\documentclass[12pt]{article}
\usepackage[top=1in,bottom=1in,left=1in,right=1in]{geometry}
\usepackage{indentfirst}
\usepackage{amsfonts,amsmath,amsthm,amssymb,hyperref,enumerate,bm}
\usepackage{longtable}
\usepackage[caption=false]{subfig}
\usepackage{leftidx}
\usepackage{lineno}
\usepackage{color,xcolor}
\usepackage{array}
\usepackage{latexsym}
\usepackage{float}
\usepackage{color}
\usepackage{longtable,supertabular,booktabs}
\usepackage{rotating}
\usepackage{cases}
\usepackage{multicol,multirow}
\usepackage{epsfig}
\usepackage{fancyhdr}       
\usepackage{dsfont}
\usepackage{hyperref}
\hypersetup{
	colorlinks=true,
	linkcolor=red,
	filecolor=blue,
	urlcolor=red,
	citecolor=blue,
}

\newtheorem{theorem}{Theorem}[section]

\newtheorem{lemma}{Lemma}[section]
\newtheorem{proposition}{Proposition}[section]
\newtheorem{corollary}{Corollary}[section]
\newtheorem{remark}{Remark}[section]

\newtheorem{example}{Example}[section]

\def\G1{G^\mathcal{C}}

 \newenvironment{prof}{\trivlist
      \item[\hskip\labelsep
      {\itshape Proof.}]\normalfont}
      {\hspace*{\fill}$\Box$\endtrivlist}
\begin{document}
\title{The matrix equation $aX^m+bY^n=cI$ over $M_2(\mathbb{Z})$}
\author{
Hongjian Li
\footnote{E-mail\,$:$ hongjian\_li@m.scnu.edu.cn}\qquad
Pingzhi Yuan
\footnote{Corresponding author. E-mail\,$:$ yuanpz@scnu.edu.cn.
Supported by National Natural Science Foundation of China (Grant No. 12171163).}\\
{\small\it  School of Mathematical Sciences, South China Normal University,}\\
{\small\it Guangzhou, Guangdong 510631, P. R. China} \\
}
\date{}
 \maketitle
\date{}

\noindent{\bf Abstract}\quad
Let $\mathbb{N}$ be the set of all positive integers and let $a,\, b,\, c$ be nonzero integers such that $\gcd\left(a,\, b,\, c\right)=1$. In this paper, we prove the following three results: (1) the solvability of the matrix equation $aX^m+bY^n=cI,\,X,\,Y\in M_2(\mathbb{Z}),\, m,\, n\in\mathbb{N}$ can be reduced to the solvability of the corresponding Diophantine equation if $XY\neq YX$ and the solvability of the equation $ax^m+by^n=c,\, m,\, n\in\mathbb{N}$ in quadratic fields if $XY=YX$; (2) we determine all non-commutative solutions of the matrix equation $X^n+Y^n=c^nI,\,X,\,Y\in M_2(\mathbb{Z}),\,n\in\mathbb{N},\,n\geq3$, and the solvability of this matrix equation can be reduced to the solvability of the equation $x^n+y^n=c^n,\, n\in\mathbb{N},\,n\geq3$ in quadratic fields if $XY=YX$; (3) we determine all solutions of the matrix equation $aX^2+bY^2=cI,\,X,\,Y\in M_2(\mathbb{Z})$.


\medskip
\noindent{\bf MR(2020) Subject Classification} 15A20, 15A24, 15B36, 11D09, 11D41

\section{Introduction}
In \cite{Vaserstein}, Vaserstein suggested solving some classical number theory problems in matrices. He considered a few classical problems of number theory with the ring $\mathbb{Z}$ substituted by the ring $M_2(\mathbb{Z})$ of $2\times2$ integral matrices, that is $2\times2$ matrices over $\mathbb{Z}$. Some classical Diophantine equations, such as Fermat's equation, Catalan's equation and Pell's equation, to matrix equations were studied by number of authors such as \cite{Chien, Cohen, Grytczuk1, Grytczuk, Khazanov, Le, Li H., Li, Qin}.

Let us recall that the Pell's equation is a Diophantine equation of the form
$$
x^2-dy^2=1,
$$
where $d$ is a positive integer which is not a perfect square. It is well-known that the Pell's equation has infinitely many solutions in positive integers $x$ and $y$. Recently, A. Grytczuk and I. Kurzyd{\l}o \cite{Grytczuk} considered the solvability of the matrix negative Pell's equation
$$
X^2-dY^2=-I,\,X,\,Y\in M_2(\mathbb{Z}),
$$
where $d$ is a positive integer. They gave a necessary and sufficient condition for the solvability of this matrix equation for nonsingular matrices $X,\,Y\in M_2(\mathbb{Z})$. In \cite{Cohen}, B. Cohen considered the solvability of the generalized matrix Pell's
equation
\begin{equation}\label{e33}
X^2-dY^2=cI,\,X,\,Y\in M_2(\mathbb{Z}),
\end{equation}
where $d$ is a square-free integer and $c$ is an arbitrary integer. He determined all solutions of equation \eqref{e33} for $c=\pm1$, as well as all non-commutative solutions for an arbitrary integer $c$. Moreover, he proposed an open problem: how about the commutative solutions of equation \eqref{e33} for an arbitrary integer $c$ ? In this paper, we give complete answers to this open problem.

The rest of this paper is organized as follows. In Section \ref{s1}, we mainly study the solvability of the matrix equation
\begin{equation}\label{e1}
aX^m+bY^n=cI,\,X,\,Y\in M_2(\mathbb{Z}),\,m,\, n\in\mathbb{N},
\end{equation}
where $a,\, b,\, c$ are nonzero integers such that $\gcd\left(a,\, b,\, c\right)=1$. Let $\lambda$ be a nonzero integer and let $n\geq3$ be a positive integer. Let $a=b=1,\,c=\lambda^n$ and $m=n$. Then equation \eqref{e1} becomes the matrix equation
\begin{equation}\label{e6}
X^n+Y^n=\lambda^nI,\,X,\,Y\in M_2(\mathbb{Z}),\,n\in\mathbb{N}, \,n\geq3.
\end{equation}
In Section \ref{s2}, we mainly study the solvability of the matrix equation \eqref{e6}. Let $m=n=2$. Then equation \eqref{e1} becomes the matrix equation
\begin{equation}\label{e7}
aX^2+bY^2=cI,\,X,\,Y\in M_2(\mathbb{Z}).
\end{equation}
In Section \ref{s3}, we mainly study the solvability of the matrix equation \eqref{e7}, and we determine all solutions of this matrix equation.

\section{The solvability of $aX^m+bY^n=cI,\,X,\,Y\in M_2(\mathbb{Z})$}\label{s1}
In this section, we will study separately commutative and non-commutative solutions of equation \eqref{e1}, i.e., solutions satisfying $XY=YX$ or $XY\neq YX$, respectively. We first study the non-commutative solutions of equation \eqref{e1}.

\begin{lemma}{\rm (\cite[Theorem 1]{Mc})}\label{le1}
Let $A=\begin{pmatrix} e & f \\ g & h \end{pmatrix}$ be an arbitrary $2\times2$ matrix and let $T=e+h$ denote its trace and $D=eh-fg$ its determinant. Let $y_n=\sum_{i=0}^{\lfloor n/2\rfloor}{n-i\choose i}T^{n-2i}(-D)^i$. Then, for $n\geq1$,
$$
A^n=\begin{pmatrix} y_n-hy_{n-1} & fy_{n-1} \\ gy_{n-1} & y_n-ey_{n-1} \end{pmatrix}.
$$
\end{lemma}

\begin{theorem}\label{p1}
Let $a,\, b,\, c$ be nonzero integers such that $\gcd\left(a,\, b,\, c\right)=1$ and let $m,\, n$ be positive integers. If there are two matrices $X,\,Y\in M_2(\mathbb{Z})$ such that
$$
aX^m+bY^n=cI,\,XY\neq YX,
$$
then $X^m$ and $Y^n$ are scalar matrices.
\end{theorem}
\begin{prof}
Note that $X^m$ is a scalar matrix if and only if $Y^n$ is a scalar matrix. So we only need to show that $Y^n$ is a scalar matrix. Let $J$ be the Jordan canonical form of $X$. Then there is a nonsingular matrix $P\in M_2(\mathbb{C})$ such that $P^{-1}XP=J$. The assumption $aX^m+bY^n=cI$ implies that $a\left(P^{-1}XP\right)^m+b\left(P^{-1}YP\right)^n=cI$, i.e.,
\begin{equation}\label{e4}
aJ^m+b\left(P^{-1}YP\right)^n=cI.
\end{equation}
Let $P^{-1}YP=\begin{pmatrix} e & f \\ g & h \end{pmatrix}$. By Lemma \ref{le1}, we have $\left(P^{-1}YP\right)^n=\begin{pmatrix} y_n-hy_{n-1} & fy_{n-1} \\ gy_{n-1} & y_n-ey_{n-1} \end{pmatrix}$, where $y_n=\sum_{i=0}^{\lfloor n/2\rfloor}{n-i\choose i}\left(tr(Y)\right)^{n-2i}\left(-\det(Y)\right)^i$. From \eqref{e4}, it follows that
\begin{equation}\label{e5}
aJ^m+b\begin{pmatrix} y_n-hy_{n-1} & fy_{n-1} \\ gy_{n-1} & y_n-ey_{n-1} \end{pmatrix}=cI.
\end{equation}

{\bf Case 1:} $J=\begin{pmatrix} x_1 & 0 \\ 0 & x_2 \end{pmatrix}$, where $x_1$ and $x_2$ are the eigenvalues of $X$.

By \eqref{e5}, we have
$$
a\begin{pmatrix} x_1^m & 0 \\ 0 & x_2^m \end{pmatrix}+b\begin{pmatrix} y_n-hy_{n-1} & fy_{n-1} \\ gy_{n-1} & y_n-ey_{n-1} \end{pmatrix}=\begin{pmatrix} c & 0 \\ 0 & c \end{pmatrix}.
$$
Comparing both sides, we have $fy_{n-1}=gy_{n-1}=0$. If $y_{n-1}\neq0$, then $f=g=0$, which implies that $P^{-1}YP=\begin{pmatrix} e & 0 \\ 0 & h \end{pmatrix}$. Since $P^{-1}XP=J=\begin{pmatrix} x_1 & 0 \\ 0 & x_2 \end{pmatrix}$, we obtain $XY=YX$, a contradiction. Therefore, $y_{n-1}=0$. Then $\left(P^{-1}YP\right)^n=y_nI$. This implies that $Y^n=y_nI$.

{\bf Case 2:} $J=\begin{pmatrix} \lambda & 1 \\ 0 & \lambda \end{pmatrix}$, where $\lambda$ is the eigenvalue of $X$.

By \eqref{e5}, we obtain
$$
a\begin{pmatrix} \lambda^m & \ast \\ 0 & \lambda^m \end{pmatrix}+b\begin{pmatrix} y_n-hy_{n-1} & fy_{n-1} \\ gy_{n-1} & y_n-ey_{n-1} \end{pmatrix}=\begin{pmatrix} c & 0 \\ 0 & c \end{pmatrix}.
$$
Comparing both sides, we have $gy_{n-1}=(e-h)y_{n-1}=0$. If $y_{n-1}\neq0$, then $g=0$ and $e=h$, which imply that $P^{-1}YP=\begin{pmatrix} e & f \\ 0 & e \end{pmatrix}$. Since $P^{-1}XP=J=\begin{pmatrix} \lambda & 1 \\ 0 & \lambda \end{pmatrix}$, we obtain $XY=YX$, a contradiction. Therefore, $y_{n-1}=0$. Then $\left(P^{-1}YP\right)^n=y_nI$. This implies that $Y^n=y_nI$.
\end{prof}

If we replace the $cI$ by $\begin{pmatrix} c_1 & 0 \\ 0 & c_2 \end{pmatrix}$ in Theorem \ref{p1} and let $XY=YX$, then we have the following proposition.

\begin{proposition}\label{p7}
Let $a,\,b$ be nonzero integers and let $m,\, n$ be positive integers. Let $c_1$ and $c_2$ be integers such that $c_1\neq c_2$. Let $X$ and $Y$ be $2\times2$ matrices over $\mathbb{Z}$. Then
\begin{equation}\label{e34}
aX^m+bY^n=\begin{pmatrix} c_1 & 0 \\ 0 & c_2 \end{pmatrix},\,XY=YX
\end{equation}
if and only if
$$
X=\begin{pmatrix} x_1 & 0 \\ 0 & x_2 \end{pmatrix},\,Y=\begin{pmatrix} y_1 & 0 \\ 0 & y_2 \end{pmatrix},
$$
where $x_1,\,x_2,\,y_1,\,y_2\in\mathbb{Z}$ satisfy $ax_i^m+by_i^n=c_i,\,i=1,\,2$.
\end{proposition}
\begin{prof}
The sufficiency is clear. We next prove necessity.

{\bf Case 1:} $X$ and $Y$ are diagonalizable.

From $XY=YX$, it follows that they are simultaneously diagonalizable. Then there is a nonsingular matrix $P=\begin{pmatrix} p_{11} & p_{12} \\ p_{21} & p_{22} \end{pmatrix}\in M_2(\mathbb{C})$ such that
$$P^{-1}XP=\begin{pmatrix} x_1 & 0 \\ 0 & x_2 \end{pmatrix},\,P^{-1}YP=\begin{pmatrix} y_1 & 0 \\ 0 & y_2 \end{pmatrix},$$
where $x_i,\,y_i,\,i=1,\,2$ are the eigenvalues of $X$ and $Y$, respectively. The assumption $aX^m+bY^n=\begin{pmatrix} c_1 & 0 \\ 0 & c_2 \end{pmatrix}$ implies that $a\left(P^{-1}XP\right)^m+b\left(P^{-1}YP\right)^n=P^{-1}\begin{pmatrix} c_1 & 0 \\ 0 & c_2 \end{pmatrix}P$. Then
$$
a\begin{pmatrix} x_1^m & 0 \\ 0 & x_2^m \end{pmatrix}+b\begin{pmatrix} y_1^n & 0 \\ 0 & y_2^n \end{pmatrix}=\frac{1}{p_{11}p_{22}-p_{12}p_{21}}\begin{pmatrix} p_{11}p_{22}c_1-p_{12}p_{21}c_2 & p_{12}p_{22}\left(c_1-c_2\right) \\ p_{11}p_{21}\left(c_2-c_1\right) & p_{11}p_{22}c_2-p_{12}p_{21}c_1 \end{pmatrix}.
$$
Comparing both sides, we have $p_{12}p_{22}=p_{11}p_{21}=0$.

If $p_{12}=0$, then $\det\left(P\right)=p_{11}p_{22}-p_{12}p_{21}=p_{11}p_{22}$. Since $\det\left(P\right)\neq0$ and $p_{11}p_{21}=0$, we have $p_{21}=0$. Then
$$
X=P\begin{pmatrix} x_1 & 0 \\ 0 & x_2 \end{pmatrix}P^{-1}=\begin{pmatrix} p_{11} & 0 \\ 0 & p_{22} \end{pmatrix}\begin{pmatrix} x_1 & 0 \\ 0 & x_2 \end{pmatrix}\begin{pmatrix} p_{11} & 0 \\ 0 & p_{22} \end{pmatrix}^{-1}=\begin{pmatrix} x_1 & 0 \\ 0 & x_2 \end{pmatrix}.
$$
Likewise, $Y=\begin{pmatrix} y_1 & 0 \\ 0 & y_2 \end{pmatrix}$. Since $X,\,Y\in M_2(\mathbb{Z})$, we have $x_1,\,x_2,\,y_1,\,y_2\in\mathbb{Z}$. The assumption $aX^m+bY^n=\begin{pmatrix} c_1 & 0 \\ 0 & c_2 \end{pmatrix}$ implies that $ax_i^m+by_i^n=c_i,\,i=1,\,2$.

If $p_{22}=0$, then $\det\left(P\right)=p_{11}p_{22}-p_{12}p_{21}=-p_{12}p_{21}$. Since $\det\left(P\right)\neq0$ and $p_{11}p_{21}=0$, we have $p_{11}=0$. Then
$$
X=P\begin{pmatrix} x_1 & 0 \\ 0 & x_2 \end{pmatrix}P^{-1}=\begin{pmatrix} 0 & p_{12} \\ p_{21} & 0 \end{pmatrix}\begin{pmatrix} x_1 & 0 \\ 0 & x_2 \end{pmatrix}\begin{pmatrix} 0 & p_{12} \\ p_{21} & 0 \end{pmatrix}^{-1}=\begin{pmatrix} x_2 & 0 \\ 0 & x_1 \end{pmatrix}.
$$
Likewise, $Y=\begin{pmatrix} y_2 & 0 \\ 0 & y_1 \end{pmatrix}$. Since $X,\,Y\in M_2(\mathbb{Z})$, we have $x_1,\,x_2,\,y_1,\,y_2\in\mathbb{Z}$. The assumption $aX^m+bY^n=\begin{pmatrix} c_1 & 0 \\ 0 & c_2 \end{pmatrix}$ implies that $ax_2^m+by_2^n=c_1$ and $ax_1^m+by_1^n=c_2$.

{\bf Case 2:} $X$ and $Y$ are not both diagonalizable.

Without loss of generality, we can assume that $X$ is not diagonalizable. Let $J$ be the Jordan canonical form of $X$. Then $J=\begin{pmatrix} \lambda & 1 \\ 0 & \lambda \end{pmatrix}$, where $\lambda$ is the eigenvalue of $X$. Moreover, there is a nonsingular matrix $P=\begin{pmatrix} p_{11} & p_{12} \\ p_{21} & p_{22} \end{pmatrix}\in M_2(\mathbb{C})$ such that $P^{-1}XP=J=\begin{pmatrix} \lambda & 1 \\ 0 & \lambda \end{pmatrix}.$ Since $XY=YX$, we have $\left(P^{-1}XP\right)\cdot\left(P^{-1}YP\right)=\left(P^{-1}YP\right)\cdot\left(P^{-1}XP\right)$, i.e., $J\cdot\left(P^{-1}YP\right)=\left(P^{-1}YP\right)\cdot J$.  This implies that $P^{-1}YP=\begin{pmatrix} y_1 & y_2 \\ 0 & y_1 \end{pmatrix}$, where $y_1,\,y_2\in\mathbb{C}$. The assumption $aX^m+bY^n=\begin{pmatrix} c_1 & 0 \\ 0 & c_2 \end{pmatrix}$ implies that $a\left(P^{-1}XP\right)^m+b\left(P^{-1}YP\right)^n=P^{-1}\begin{pmatrix} c_1 & 0 \\ 0 & c_2 \end{pmatrix}P$. Then
$$
a\begin{pmatrix} \lambda^m & \ast \\ 0 & \lambda^m \end{pmatrix}+b\begin{pmatrix} y_1^n & \ast \\ 0 & y_1^n \end{pmatrix}=\frac{1}{p_{11}p_{22}-p_{12}p_{21}}\begin{pmatrix} p_{11}p_{22}c_1-p_{12}p_{21}c_2 & p_{12}p_{22}\left(c_1-c_2\right) \\ p_{11}p_{21}\left(c_2-c_1\right) & p_{11}p_{22}c_2-p_{12}p_{21}c_1 \end{pmatrix}.
$$
Comparing both sides, we have $p_{11}p_{21}=p_{11}p_{22}+p_{12}p_{21}=0$. This implies that $\det\left(P\right)=p_{11}p_{22}-p_{12}p_{21}=0$, a contradiction.
\end{prof}

\begin{remark}
\rm All commutative solutions of equation \eqref{e34} are given by Proposition \ref{p7}. How about the non-commutative solutions of equation \eqref{e34}? This is an interesting problem that lies out of the scope of this paper.
\end{remark}

About scalar matrices, we have the following lemma and proposition.

\begin{lemma}{\rm (\cite{Zhong})}\label{le2}
Let $X$ be a $2\times2$ matrix over $\mathbb{Z}$ such that $X^n$ is a scalar matrix for some $n\in\mathbb{N}$, and let $k$ be the smallest positive integer with such property. Then the following statements hold.
\begin{enumerate}
\item[\rm1)]\, $k\in\{1,\,2,\,3,\,4,\,6\}$;
\item[\rm2)]\,
\begin{enumerate}
\item[\rm(i)]\, $k=1$ if and only if $X=aI,\,a\in\mathbb{Z}$;
\item[\rm(ii)]\,$k=2$ if and only if $X=\begin{pmatrix} a & b \\ c & -a \end{pmatrix}$, where $a,\, b,\, c\in\mathbb{Z}$ satisfy $a^2+b^2+c^2\neq0$. Moreover, $X^2=\left(a^2+bc\right)I$;
\item[\rm(iii)]\,$k=3$ if and only if $X=\begin{pmatrix} a & b \\ c & d \end{pmatrix}$, where $a,\, b,\, c,\,d\in\mathbb{Z}$ satisfy $\left(a+d\right)^2=ad-bc$ and $a+d\neq0$. Moreover, $X^3=-\left(a+d\right)^3I$;
\item[\rm(iv)]\,$k=4$ if and only if $X=\begin{pmatrix} a & b \\ c & d \end{pmatrix}$, where $a,\, b,\, c,\,d\in\mathbb{Z}$ satisfy  $\left(a+d\right)^2=2\left(ad-bc\right)$ and $a+d\neq0$. Moreover, $X^4=-4\left(\frac{a+d}{2}\right)^4I=-\left(ad-bc\right)^2I$;
\item[\rm(v)]\,$k=6$ if and only if $X=\begin{pmatrix} a & b \\ c & d \end{pmatrix}$, where $a,\, b,\, c,\,d\in\mathbb{Z}$ satisfy  $\left(a+d\right)^2=3\left(ad-bc\right)$ and $a+d\neq0$. Moreover, $X^6=-\left(ad-bc\right)^3I$.
\end{enumerate}
\end{enumerate}
\end{lemma}

\begin{proposition}\label{p3}
Let $X\in M_2(\mathbb{Z})$ be a nonsingular matrix such that $X^n$ is a scalar matrix for some $n\in\mathbb{N}$, and let $k$ be the smallest positive integer with such property. Then for $m\in\mathbb{N}$, $X^m$ is a scalar matrix if and only if $k\mid m$.
\end{proposition}
\begin{prof}
The sufficiency is clear. We next prove necessity. Since $X^m$ and $X^k$ are scalar matrices, we have $X^m=\lambda I$ and $X^k=\mu I$ for some $\lambda,\,\mu\in\mathbb{Z}\backslash\{0\}$. Let $m=kq+r$, where $q,\, r\in\mathbb{Z}$ and $0\leq r<k$. Then
$$
X^r=X^{m-kq}=X^m\cdot \left(X^k\right)^{-q}=\lambda I\cdot\left(\mu I\right)^{-q}=\frac{\lambda}{\mu^q}I.
$$
If $r\neq0$, then we obtain a contradiction to the minimality of $k$. Thus, $r=0$. This means that $k\mid m$.
\end{prof}

About commutative $2\times2$ integral matrices, we have the following lemma.

\begin{lemma}\label{le3}
Let $X=\begin{pmatrix} t_1 & t_2 \\ t_3 & t_4 \end{pmatrix}$ and $Y=\begin{pmatrix} s_1 & s_2 \\ s_3 & s_4 \end{pmatrix}$ be $2\times2$ matrices over $\mathbb{Z}$. Then $XY=YX$ if and only if the vectors $\overrightarrow{t}=(t_1-t_4,\,t_2,\,t_3)$ and $\overrightarrow{s}=(s_1-s_4,\,s_2,\,s_3)$ are linearly dependent over $\mathbb{Q}$.
\end{lemma}
\begin{prof}
By a direct computation, we have
$$
\begin{aligned}
XY-YX
&=\begin{pmatrix} t_2s_3-s_2t_3 & (t_1-t_4)s_2-(s_1-s_4)t_2 \\ (s_1-s_4)t_3-(t_1-t_4)s_3 & s_2t_3-t_2s_3 \end{pmatrix}\\
&=\begin{pmatrix}
\begin{vmatrix}
t_2 & t_3  \\
s_2 & s_3  \\
\end{vmatrix} &
\begin{vmatrix}
t_1-t_4 & t_2  \\
s_1-s_4 & s_2  \\
\end{vmatrix} \\
-\begin{vmatrix}
t_1-t_4 & t_3  \\
s_1-s_4 & s_3  \\
\end{vmatrix} &
-\begin{vmatrix}
t_2 & t_3  \\
s_2 & s_3  \\
\end{vmatrix}
\end{pmatrix}.
\end{aligned}
$$
Let $\overrightarrow{i}=\left(1,\,0,\,0\right)$, $\overrightarrow{j}=\left(0,\,1,\,0\right)$ and $\overrightarrow{k}=\left(0,\,0,\,1\right)$. Then $XY=YX$ if and only if
$$\begin{vmatrix}
t_2 & t_3  \\
s_2 & s_3  \\
\end{vmatrix}=\begin{vmatrix}
t_1-t_4 & t_2  \\
s_1-s_4 & s_2  \\
\end{vmatrix}=\begin{vmatrix}
t_1-t_4 & t_3  \\
s_1-s_4 & s_3  \\
\end{vmatrix}=0$$
if and only if
$$
\overrightarrow{t}\times\overrightarrow{s}=
\begin{vmatrix}
\overrightarrow{i}&\overrightarrow{j}&\overrightarrow{k}\\
t_1-t_4 & t_2 & t_3\\
s_1-s_4 & s_2 & s_3\\
\end{vmatrix}=\begin{vmatrix}
t_2 & t_3  \\
s_2 & s_3  \\
\end{vmatrix}\overrightarrow{i}-\begin{vmatrix}
t_1-t_4 & t_3  \\
s_1-s_4 & s_3  \\
\end{vmatrix}\overrightarrow{j}+\begin{vmatrix}
t_1-t_4 & t_2  \\
s_1-s_4 & s_2  \\
\end{vmatrix}\overrightarrow{k}=\overrightarrow{0}
$$
if and only if the vectors $\overrightarrow{t}$ and $\overrightarrow{s}$ are linearly dependent over $\mathbb{Q}$.
\end{prof}

From Theorem \ref{p1}, Lemmas \ref{le2}, \ref{le3} and Proposition \ref{p3}, we conclude that finding the non-commutative solutions of equation \eqref{e1} can be reduced to finding the solutions of the corresponding Diophantine equation. Next, we give an example to illustrate it, i.e., Proposition \ref{p2}.

\begin{proposition}\label{p2}
Let $a,\, b,\, c$ be nonzero integers such that $\gcd\left(a,\, b,\, c\right)=1$. Let $X$ and $Y$ be $2\times2$ matrices over $\mathbb{Z}$. Then the following statements hold.
\begin{enumerate}
\item[\rm1)]\,
$$
aX^2+bY^2=cI,\,XY\neq YX
$$
if and only if
$$
X=\begin{pmatrix} t_1 & t_2 \\ t_3 & -t_1 \end{pmatrix},\,Y=\begin{pmatrix} s_1 & s_2 \\ s_3 & -s_1 \end{pmatrix},
$$
where $t_1,\,t_2,\,t_3,\,s_1,\,s_2,\,s_3\in\mathbb{Z}$ satisfy $a\left(t_1^2+t_2t_3\right)+b\left(s_1^2+s_2s_3\right)=c$ and the vectors
$\overrightarrow{t}=(t_1,\,t_2,\,t_3)$ and $\overrightarrow{s}=(s_1,\,s_2,\,s_3)$ are linearly independent over $\mathbb{Q}$;
\item[\rm2)]\,
$$
X^4+Y^4=c^4I,\,XY\neq YX
$$
if and only if
$$
X=\begin{pmatrix} t_1 & t_2 \\ t_3 & -t_1 \end{pmatrix},\,Y=\begin{pmatrix} s_1 & s_2 \\ s_3 & -s_1 \end{pmatrix},
$$
where $t_1,\,t_2,\,t_3,\,s_1,\,s_2,\,s_3\in\mathbb{Z}$ satisfy $\left(t_1^2+t_2t_3\right)^2+\left(s_1^2+s_2s_3\right)^2=c^4$ and the vectors $\overrightarrow{t}=(t_1,\,t_2,\,t_3)$ and $\overrightarrow{s}=(s_1,\,s_2,\,s_3)$ are linearly independent over $\mathbb{Q}$.
\end{enumerate}
\end{proposition}

\begin{prof}
1) Sufficiency follows from Lemma \ref{le2} 2) (ii) and Lemma \ref{le3}. We next prove necessity. From Theorem \ref{p1}, it follows that $X^2$ and $Y^2$ are scalar matrices. Since $XY\neq YX$, $X$ and $Y$ are not scalar matrices. By Lemma \ref{le2} 2) (ii), we have
$$
X=\begin{pmatrix} t_1 & t_2 \\ t_3 & -t_1 \end{pmatrix},\,Y=\begin{pmatrix} s_1 & s_2 \\ s_3 & -s_1 \end{pmatrix},
$$
where $t_1,\,t_2,\,t_3,\,s_1,\,s_2,\,s_3$ are integers. Moreover, $X^2=\left(t_1^2+t_2t_3\right)I$ and $Y^2=\left(s_1^2+s_2s_3\right)I$. The assumption $aX^2+bY^2=cI$ implies that $a\left(t_1^2+t_2t_3\right)+b\left(s_1^2+s_2s_3\right)=c$. Let $\overrightarrow{t}=(t_1,\,t_2,\,t_3)$ and $\overrightarrow{s}=(s_1,\,s_2,\,s_3)$. Then it follows from Lemma \ref{le3} that $\overrightarrow{t}$ and $\overrightarrow{s}$ are linearly independent over $\mathbb{Q}$.

2) Sufficiency follows from Lemma \ref{le2} 2) (ii) and Lemma \ref{le3}. We next prove necessity. From Theorem \ref{p1}, it follows that $X^4$ and $Y^4$ are scalar matrices. Let $k$ and $l$ be the smallest positive integers such that $X^k$ and $Y^l$ are scalar matrices, respectively. Since $XY\neq YX$, we have $k,\,l\neq1$. By Lemma \ref{le2} and Proposition \ref{p3}, we have $k,\,l\in\{2,\,4\}$.

{\bf Case 1:} $k=l=2$.

By Lemma \ref{le2}, we have
$$
X=\begin{pmatrix} t_1 & t_2 \\ t_3 & -t_1 \end{pmatrix},\,Y=\begin{pmatrix} s_1 & s_2 \\ s_3 & -s_1 \end{pmatrix},
$$
where $t_1,\,t_2,\,t_3,\,s_1,\,s_2,\,s_3$ are integers. Moreover, $X^2=\left(t_1^2+t_2t_3\right)I$ and $Y^2=\left(s_1^2+s_2s_3\right)I$. The assumption $X^4+Y^4=c^4I$ implies that $\left(t_1^2+t_2t_3\right)^2+\left(s_1^2+s_2s_3\right)^2=c^4$. Let $\overrightarrow{t}=(t_1,\,t_2,\,t_3)$ and $\overrightarrow{s}=(s_1,\,s_2,\,s_3)$. Then it follows from Lemma \ref{le3} that $\overrightarrow{t}$ and $\overrightarrow{s}$ are linearly independent over $\mathbb{Q}$.

{\bf Case 2:} $k=2,\,l=4$.

By Lemma \ref{le2}, we have $X^2=uI$ and $Y^4=-4v^4I$ for some $u\in\mathbb{Z}$ and $v\in\mathbb{Z}\backslash\{0\}$. The assumption $X^4+Y^4=c^4I$ implies that $c^4+4v^4=u^2$, which is impossible.

{\bf Case 3:} $k=4,\,l=2$.

This case is also impossible, where the reason is similar to Case 2.

{\bf Case 4:} $k=l=4$.

By Lemma \ref{le2}, we have $X^4=-4u^4I$ and $Y^4=-4v^4I$ for some $u,\,v\in\mathbb{Z}\backslash\{0\}$. The assumption $X^4+Y^4=c^4I$ implies that $-4u^4-4v^4=c^4$, which is impossible.
\end{prof}

We now study the commutative solutions of equation \eqref{e1}.
\begin{proposition}\label{p8}
If $\left(X,\,Y\right)$ is a solution of equation \eqref{e1}, then
\begin{equation}\label{e35}
ax_i^m+by_i^n=c,\,i=1,\,2,
\end{equation}
where $x_i,\,y_i,\,i=1,\,2$ are the eigenvalues of $X$ and $Y$, respectively.
\end{proposition}
\begin{prof}
Suppose that $\begin{pmatrix} x_1 & \ast \\ 0 & x_2 \end{pmatrix}$ is the Jordan canonical form of $X$. Then there is a nonsingular matrix $P\in M_2(\mathbb{C})$ such that $P^{-1}XP=\begin{pmatrix} x_1 & \ast \\ 0 & x_2 \end{pmatrix}$. The assumption $aX^m+bY^n=cI$ implies that $a\left(P^{-1}XP\right)^m+b\left(P^{-1}YP\right)^n=cI$. Then we obtain $\begin{pmatrix} ax_1^m & \ast \\ 0 & ax_2^m \end{pmatrix}+b\left(P^{-1}YP\right)^n=cI$, which implies that
\begin{equation}\label{e41}
bY^n=P\begin{pmatrix} c-ax_1^m & \ast \\ 0 & c-ax_2^m \end{pmatrix}P^{-1}.
\end{equation}
Comparing the eigenvalues of both sides of \eqref{e41}, we have $ax_i^m+by_i^n=c,\,i=1,\,2$.
\end{prof}

If the eigenvalues of $X$ and $Y$ are integers, then equation \eqref{e35} becomes the Diophantine equation $ax^m+by^n=c,\,x,\,y\in\mathbb{Z}$. So we can assume that the eigenvalues of $X$ or $Y$ are not integers. Let $A=\begin{pmatrix} e & f \\ g & 0 \end{pmatrix}\in M_2\left(\mathbb{Z}\right)$ be a given matrix such that $fg\neq0$ and $\gcd(e,\, f,\, g)=1$, and let $C(A)=\{B\in M_2(\mathbb{Z}): AB=BA\}$. In \cite{Li}, we showed that the solvability of the matrix equation
\begin{equation}\label{e36}
aX^m+bY^n=cI,\,XY=YX,\, m,\, n\in\mathbb{N}
\end{equation}
in $M_2(\mathbb{Z})$ can be reduced to the solvability of the matrix equation
\begin{equation}\label{e2}
aX^m+bY^n=cI,\,m,\, n\in\mathbb{N}
\end{equation}
in $C(A)$, and finally to the solvability of the equation
\begin{equation}\label{e3}
ax^m+by^n=c,\, m,\, n\in\mathbb{N}
\end{equation}
in quadratic fields. As a corollary of \cite[Theorem 3.1]{Li}, we have the following theorem.

\begin{theorem}\label{th1}
Let $A=\begin{pmatrix} e & f \\ g & 0 \end{pmatrix}\in M_2\left(\mathbb{Z}\right)$ be a given matrix such that $fg\neq0$ and $\gcd(e,\, f,\, g)=1$. Let $K=\mathbb{Q}\left(\sqrt{e^2+4fg}\right)$ and let $\mathcal{O}_K$ be its ring of integers. Then the following statements hold.
\begin{enumerate}
\item[\rm1)]\, If $e^2+4fg$ is a square, then equation \eqref{e2}  has a non-trivial solution in $C(A)$ if and only if equation \eqref{e3} has a non-trivial solution in $\mathbb{Z}$;
\item[\rm2)]\, If $e^2+4fg$ is not a square and $D$ is the unique square-free integer such that $e^2+4fg=k^2D$ for some $k\in\mathbb{N}$, then equation \eqref{e2} has a non-trivial solution in $C(A)$ if and only if equation \eqref{e3} has a non-trivial solution $\left(x,\, y\right)$ in $\mathcal{O}_K$ such that $x,\, y$ can be written in the form
$$
\frac{s+t\sqrt{D}}{2},\quad s,\, t\in\mathbb{Z},\, k\mid t.
$$
\end{enumerate}
\end{theorem}

From Theorem \ref{th1}, we conclude that the solvability of equation \eqref{e36} in $M_2(\mathbb{Z})$ can be reduced to the solvability of equation \eqref{e3} in quadratic fields. However, the solvability of equation \eqref{e3} in quadratic fields is unsolved.

\section{The solvability of $X^n+Y^n=\lambda^nI,\,X,\,Y\in M_2(\mathbb{Z})$}\label{s2}

In this section, we mainly consider the non-trivial solutions of equation \eqref{e6}, i.e., solutions satisfying $\det\left(XY\right)\neq0$. Indeed, suppose that $\left(X,\,Y\right)$ is a solution of equation \eqref{e6} such that $\det\left(XY\right)=0$. Without loss of generality, we can assume that $\det\left(X\right)=0$. If the eigenvalues of $X$ are both equal to zero, then $X^2=O$. In this case, equation \eqref{e6} becomes the matrix equation
\[
Y^n=\lambda^nI,\,Y\in M_2(\mathbb{Z}),\,n\in\mathbb{N},\,n\geq3.
\]
However, all solutions of this matrix equation can be given by Lemma \ref{le2}. So we only need to consider the non-trivial solutions of equation \eqref{e6}. In \cite{Li H.}, we showed that equation \eqref{e6} has no non-trivial solutions if the eigenvalues of $X$ or $Y$ are integers. Moreover, we gave all non-trivial solutions of equation \eqref{e6} for $\lambda=\pm1$. In this section, we will study separately commutative and non-commutative solutions of equation \eqref{e6} for an arbitrary nonzero integer $\lambda$. We first study the non-commutative solutions of equation \eqref{e6}.

\begin{theorem}\label{th4}
Equation \eqref{e6} has no non-commutative non-trivial solutions for $n\neq4$.
\end{theorem}
\begin{prof}
Suppose that $\left(X,\,Y\right)$ is a non-commutative non-trivial solution of equation \eqref{e6} for $n\neq4$. Then we have $X^n+Y^n=\lambda^nI,\,XY\neq YX$ and $\det\left(XY\right)\neq0$. From Theorem \ref{p1}, it follows that $X^n$ and $Y^n$ are scalar matrices. Let $k$ and $l$ be the smallest positive integers such that $X^k$ and $Y^l$ are scalar matrices, respectively. Since $XY\neq YX$, we have $k,\,l\neq1$. By Lemma \ref{le2}, we have $k,\,l\in\{2,\,3,\,4,\,6\}$. From Proposition \ref{p3}, it follows that $k\mid n$ and $l\mid n$. Let $q$ be the least common multiple of $k$ and $l$. Then $q\mid n$. Since $k,\,l\in\{2,\,3,\,4,\,6\}$, we have $q\in\{2,\,3,\,4,\,6,\,12\}$. If $\gcd\left(n,\,6\right)=1$, then $q\nmid n$, a contradiction.

{\bf Case 1:} $n\equiv0\pmod6$. Then $n=6m$ for some positive integer $m$.

{\bf Subcase 1.1:} $2\mid m$.

Then $m=2t$ for some positive integer $t$, which implies that $n=12t$. By Lemma \ref{le2}, we have $X^{12}=a^3I$ and $Y^{12}=b^3I$ for some $a,\,b\in\mathbb{Z}\backslash\{0\}$. The assumption $X^n+Y^n=\lambda^nI$ implies that $a^{3t}I+b^{3t}I=\lambda^{12t}I$. Then $\left(a^t\right)^3+\left(b^t\right)^3=\left(\lambda^{4t}\right)^3$, which is impossible by Fermat's last theorem.

{\bf Subcase 1.2:} $2\nmid m$.

By Proposition \ref{p3}, we have $k,\,l\in\{2,\,3,\,6\}$. By Lemma \ref{le2}, we obtain $X^{6}=a^3I$ and $Y^{6}=b^3I$ for some $a,\,b\in\mathbb{Z}\backslash\{0\}$. The assumption $X^n+Y^n=\lambda^nI$ implies that $a^{3m}I+b^{3m}I=\lambda^{6m}I$. Then $\left(a^m\right)^3+\left(b^m\right)^3=\left(\lambda^{2m}\right)^3$, which is impossible by Fermat's last theorem.

{\bf Case 2:} $n\equiv2\pmod6$. Then $n=2+6m$ for some positive integer $m$.

{\bf Subcase 2.1:} $2\mid m$.

Then $m=2t$ for some positive integer $t$, which implies that $n=2+12t$. By Proposition \ref{p3}, we have $k=l=2$. By Lemma \ref{le2}, we have $X^{2}=aI$ and $Y^{2}=bI$ for some $a,\,b\in\mathbb{Z}\backslash\{0\}$. The assumption $X^n+Y^n=\lambda^nI$ implies that $a^{1+6t}I+b^{1+6t}I=\lambda^{2+12t}I$. Then $a^{1+6t}+b^{1+6t}=\left(\lambda^2\right)^{1+6t}$, which is impossible by Fermat's last theorem.

{\bf Subcase 2.2:} $2\nmid m$.

Then $1+3m=2t$ for some positive integer $t\geq2$. Moreover, $3\nmid t$. We obtain $n=2+6m=2\left(1+3m\right)=4t$. By Proposition \ref{p3}, we have $k,\,l\in\{2,\,4\}$. By Lemma \ref{le2}, we obtain $X^{4}=\pm a^2I$ and $Y^{4}=\pm b^2I$ for some $a,\,b\in\mathbb{Z}\backslash\{0\}$. The assumption $X^n+Y^n=\lambda^nI$ implies that $\left(\pm a^2\right)^tI+\left(\pm b^2\right)^tI=\left(\lambda^{4}\right)^tI$. Then $\left(\pm a^2\right)^t+\left(\pm b^2\right)^t=\left(\lambda^{4}\right)^t$, which is impossible by Fermat's last theorem.

{\bf Case 3:} $n\equiv3\pmod6$.

Then $n=3+6m$ for some non-negative integer $m$.  By Proposition \ref{p3}, we have $k=l=3$. By Lemma \ref{le2}, we have $X^{3}=a^3I$ and $Y^{3}=b^3I$ for some $a,\,b\in\mathbb{Z}\backslash\{0\}$. The assumption $X^n+Y^n=\lambda^nI$ implies that $a^{n}I+b^{n}I=\lambda^{n}I$. Then $a^{n}+b^{n}=\lambda^{n}$, which is impossible by Fermat's last theorem.

{\bf Case 4:} $n\equiv4\pmod6$. Since $n\neq4$, we have $n=4+6m$ for some positive integer $m$.

{\bf Subcase 4.1:} $2\nmid m$.

By Proposition \ref{p3}, we have $k=l=2$. By Lemma \ref{le2}, we obtain $X^{2}=aI$ and $Y^{2}=bI$ for some $a,\,b\in\mathbb{Z}\backslash\{0\}$. The assumption $X^n+Y^n=\lambda^nI$ implies that $a^{2+3m}I+b^{2+3m}I=\left(\lambda^{2}\right)^{2+3m}I$. Then $a^{2+3m}+b^{2+3m}=\left(\lambda^{2}\right)^{2+3m}$, which is impossible by Fermat's last theorem.

{\bf Subcase 4.2:} $2\mid m$.

Then $m=2t$ for some positive integer $t$, which implies that $n=4+12t$. By Proposition \ref{p3}, we have $k,\,l\in\{2,\,4\}$. By Lemma \ref{le2}, we have $X^{4}=aI$ and $Y^{4}=bI$ for some $a,\,b\in\mathbb{Z}\backslash\{0\}$. The assumption $X^n+Y^n=\lambda^nI$ implies that $a^{1+3t}I+b^{1+3t}I=\lambda^{4+12t}I$. Then $a^{1+3t}+b^{1+3t}=\left(\lambda^4\right)^{1+3t}$, which is impossible by Fermat's last theorem.
\end{prof}

\begin{remark}
\rm All non-commutative solutions of equation \eqref{e6} are given by Proposition \ref{p2} 2) for $n=4$.
\end{remark}

We now study the commutative solutions of equation \eqref{e6}. If the eigenvalues of $X$ and $Y$ are integers, then it follows from Proposition \ref{p8} that
\[
x_i^n+y_i^n=\lambda^n,\,i=1,\,2,\,n\in\mathbb{N}, \,n\geq3,
\]
where $x_i,\,y_i\in\mathbb{Z}\backslash\{0\},\,i=1,\,2$ are the eigenvalues of $X$ and $Y$, respectively. We know that this is impossible by Fermat's last theorem. Thus, we can assume that the eigenvalues of $X$ or $Y$ are not integers. Let $A=\begin{pmatrix} e & f \\ g & 0 \end{pmatrix}\in M_2\left(\mathbb{Z}\right)$ be a given matrix such that $fg\neq0$ and $\gcd(e,\, f,\, g)=1$, and let $C(A)=\{B\in M_2(\mathbb{Z}): AB=BA\}$. In \cite{Li}, we showed that the solvability of the matrix equation
\begin{equation}\label{e37}
X^n+Y^n=\lambda^nI,\,XY=YX,\,n\in\mathbb{N}, \,n\geq3
\end{equation}
in $M_2(\mathbb{Z})$ can be reduced to the solvability of the matrix equation
\begin{equation}\label{e8}
X^n+Y^n=\lambda^nI,\, n\in\mathbb{N}, \,n\geq3
\end{equation}
in $C(A)$, and finally to the solvability of the equation
\begin{equation}\label{e9}
x^n+y^n=\lambda^n,\, n\in\mathbb{N},\,n\geq3
\end{equation}
in quadratic fields. As a corollary of \cite[Theorem 3.1]{Li}, we have the following theorem.

\begin{theorem}\label{th2}
Let $A=\begin{pmatrix} e & f \\ g & 0 \end{pmatrix}\in M_2\left(\mathbb{Z}\right)$ be a given matrix such that $fg\neq0$ and $\gcd(e,\, f,\, g)=1$. Let $K=\mathbb{Q}\left(\sqrt{e^2+4fg}\right)$ and let $\mathcal{O}_K$ be its ring of integers. Then the following statements hold.
\begin{enumerate}
\item[\rm1)]\, If $e^2+4fg$ is a square, then equation \eqref{e8} has no non-trivial solutions in $C(A)$;
\item[\rm2)]\, If $e^2+4fg$ is not a square and $D$ is the unique square-free integer such that $e^2+4fg=k^2D$ for some $k\in\mathbb{N}$, then equation \eqref{e8} has a non-trivial solution in $C(A)$ if and only if equation \eqref{e9} has a non-trivial solution $\left(x,\, y\right)$ in $\mathcal{O}_K$ such that $x,\, y$ can be written in the form
$$
\frac{s+t\sqrt{D}}{2},\quad s,\, t\in\mathbb{Z},\, k\mid t.
$$
\end{enumerate}
\end{theorem}

From Theorem \ref{th2}, we conclude that the solvability of equation \eqref{e37} in $M_2(\mathbb{Z})$ can be reduced to the solvability of equation \eqref{e9} in quadratic fields. However, the solvability of equation \eqref{e9} in quadratic fields is unsolved. We next list a known result about the solvability of the Fermat's equation in quadratic fields.

\begin{lemma}{\rm (\cite{Aigner})}\label{le4}
The Fermat's equation
$$
x^n+y^n=z^n,\, n\in\mathbb{N}, \,n\geq3
$$
has no non-trivial solutions in quadratic fields for $n=6,\, 9$.
\end{lemma}

By Theorems \ref{th4}, \ref{th2} and Lemma \ref{le4}, we have the following proposition.

\begin{proposition}\label{p9}
Equation \eqref{e6} has no non-trivial solutions in $M_2(\mathbb{Z})$ for $n=6,\, 9$.
\end{proposition}

\begin{prof}
Directly from Theorems \ref{th4}, \ref{th2} and Lemma \ref{le4}.
\end{prof}

\begin{corollary}
Let $\lambda$ be a nonzero integer and let $i,\,j,\,k$ be positive integers such that $6\mid\gcd\left(i,\,j,\,k\right)$ or $9\mid\gcd\left(i,\,j,\,k\right)$. Then the matrix equation
\[
X^i+Y^j=\lambda^kI
\]
has no non-trivial solutions in $M_2(\mathbb{Z})$.
\end{corollary}

\begin{prof}
Directly from Proposition \ref{p9}.
\end{prof}

\section{All solutions of $aX^2+bY^2=cI,\,X,\,Y\in M_2(\mathbb{Z})$}\label{s3}
Note that all non-commutative solutions of equation \eqref{e7} are given by Proposition \ref{p2} 1). Thus, in this section, we only need to study the commutative solutions of equation \eqref{e7}.

\begin{theorem}\label{th3}
Let $a,\, b,\, c$ be nonzero integers such that $-ab$ is not a square and $\gcd\left(a,\, b,\, c\right)=1$. Let $X$ and $Y$ be $2\times2$ matrices over $\mathbb{Z}$. Then
$$
aX^2+bY^2=cI,\,XY=YX
$$
if and only if one of the following statements holds.
\begin{enumerate}
\item[\rm(i)]\, $X=t_1I$ and $Y=t_2I$, where $t_1,\,t_2\in\mathbb{Z}$ satisfy $at_1^2+bt_2^2=c$;
\item[\rm(ii)]\, $X=t_1I$ and $Y=\begin{pmatrix} t_4 & t_2 \\ t_3 & -t_4 \end{pmatrix}$, where $t_1,\,t_2,\,t_3,\,t_4\in\mathbb{Z}$ satisfy $at_1^2+b\left(t_4^2+t_2t_3\right)=c$;
\item[\rm(iii)]\,$X=t_1I+\frac{u-c}{g}\begin{pmatrix} 0 & t_2 \\ t_3 & -k \end{pmatrix}$ and $Y=t_4I+\frac{va}{g}\begin{pmatrix} 0 & t_2 \\ t_3 & -k \end{pmatrix}$, where $t_1,\,t_2,\,t_3,\,t_4,\,u,\,v,\,k\in\mathbb{Z},\,u\neq c,\,g=\gcd\left(va,\,u-c\right)$ and the following relations hold:
$$
u^2+v^2ab=c^2,\,at_1^2+bt_4^2+\frac{2act_2t_3}{g^2}\left(c-u\right)=c,\,\left(gt_1+ck\right)\left(u-c\right)+vbgt_4=0.
$$
\end{enumerate}
\end{theorem}

\begin{prof}
We now prove sufficiency. We need to verify the following three cases.

First, suppose that $X=t_1I$ and $Y=t_2I$, where $t_1,\,t_2\in\mathbb{Z}$ satisfy $at_1^2+bt_2^2=c$. Then $XY=YX$ and $aX^2+bY^2=at_1^2I+bt_2^2I=cI$.

Then, suppose that $X=t_1I$ and $Y=\begin{pmatrix} t_4 & t_2 \\ t_3 & -t_4 \end{pmatrix}$, where $t_1,\,t_2,\,t_3,\,t_4\in\mathbb{Z}$ satisfy $at_1^2+b\left(t_4^2+t_2t_3\right)=c$. Then $XY=YX$ and $Y^2=\left(t_4^2+t_2t_3\right)I$. Moreover,
$$aX^2+bY^2=at_1^2I+b\left(t_4^2+t_2t_3\right)I=\left(at_1^2+b\left(t_4^2+t_2t_3\right)\right)I=cI.$$

Finally, suppose that $X=t_1I+\frac{u-c}{g}\begin{pmatrix} 0 & t_2 \\ t_3 & -k \end{pmatrix}$ and $Y=t_4I+\frac{va}{g}\begin{pmatrix} 0 & t_2 \\ t_3 & -k \end{pmatrix}$, where $t_1,\,t_2,\,t_3,\,t_4,\,u,\,v,\,k\in\mathbb{Z}$ satisfy the above conditions. Then $XY=YX$. Moreover,
$$
\begin{aligned}
aX^2+bY^2=&\left(at_1^2+bt_4^2+\frac{2act_2t_3}{g^2}\left(c-u\right)\right)I+\frac{2a}{g^2}\left[\left(gt_1+ck\right)\left(u-c\right)+vbgt_4\right]\begin{pmatrix} 0 & t_2 \\ t_3 & -k \end{pmatrix}\\
=&\left(at_1^2+bt_4^2+\frac{2act_2t_3}{g^2}\left(c-u\right)\right)I=cI.
\end{aligned}
$$

We next prove necessity. Let $X=\begin{pmatrix} x_1 & x_2 \\ x_3 & x_4 \end{pmatrix}$ and $Y=\begin{pmatrix} y_1 & y_2 \\ y_3 & y_4 \end{pmatrix}$. Then the assumption $aX^2+bY^2=cI$ implies that $a^2X^2+abY^2=acI$. Since $XY=YX$, we have
\begin{equation}\label{e10}
\left(aX+\sqrt{-ab}Y\right)\left(aX-\sqrt{-ab}Y\right)=acI.
\end{equation}
Note that
$$
aX+\sqrt{-ab}Y=\begin{pmatrix} ax_1+y_1\sqrt{-ab} & ax_2+y_2\sqrt{-ab} \\ ax_3+y_3\sqrt{-ab} & ax_4+y_4\sqrt{-ab} \end{pmatrix}
$$
and
$$
aX-\sqrt{-ab}Y=\begin{pmatrix} ax_1-y_1\sqrt{-ab} & ax_2-y_2\sqrt{-ab} \\ ax_3-y_3\sqrt{-ab} & ax_4-y_4\sqrt{-ab} \end{pmatrix}.
$$
By a direct computation, we have
\[
\det\left(aX+\sqrt{-ab}Y\right)=a\left(u+v\sqrt{-ab}\right)\,\mbox{~and~}\,\det\left(aX-\sqrt{-ab}Y\right)=a\left(u-v\sqrt{-ab}\right),
\]
where $u=a\left(x_1x_4-x_2x_3\right)-b\left(y_1y_4-y_2y_3\right)$ and $v=x_1y_4+x_4y_1-x_2y_3-x_3y_2$. By computing the determinants of both sides of  \eqref{e10}, we have $u^2+v^2ab=c^2$. From \eqref{e10}, it follows that
\begin{equation}\label{e11}
\left(aX+\sqrt{-ab}Y\right)^{-1}=\frac{1}{ac}\left(aX-\sqrt{-ab}Y\right)=\frac{1}{ac}\begin{pmatrix} ax_1-y_1\sqrt{-ab} & ax_2-y_2\sqrt{-ab} \\ ax_3-y_3\sqrt{-ab} & ax_4-y_4\sqrt{-ab} \end{pmatrix}.
\end{equation}
Moreover,
\begin{equation}\label{e12}
\begin{aligned}
\left(aX+\sqrt{-ab}Y\right)^{-1}&=\frac{1}{\det\left(aX+\sqrt{-ab}Y\right)}{\rm adj}\left(aX+\sqrt{-ab}Y\right)\\
&=\frac{1}{a\left(u+v\sqrt{-ab}\right)}\begin{pmatrix} ax_4+y_4\sqrt{-ab} & -\left(ax_2+y_2\sqrt{-ab}\right) \\ -\left(ax_3+y_3\sqrt{-ab}\right) & ax_1+y_1\sqrt{-ab} \end{pmatrix},
\end{aligned}
\end{equation}
where ${\rm adj}\left(aX+\sqrt{-ab}Y\right)$ is the adjugate of $aX+\sqrt{-ab}Y$. Comparing \eqref{e11} and \eqref{e12}, we have
\begin{subnumcases}
{}acx_1+\left(-y_1c\right)\sqrt{-ab}=\left(aux_4+vaby_4\right)+\left(uy_4-avx_4\right)\sqrt{-ab}, \label{a}\\
\left(-acx_2\right)+y_2c\sqrt{-ab}=\left(aux_2+vaby_2\right)+\left(uy_2-avx_2\right)\sqrt{-ab}, \label{b} \\
\left(-acx_3\right)+y_3c\sqrt{-ab}=\left(aux_3+vaby_3\right)+\left(uy_3-avx_3\right)\sqrt{-ab}, \label{c} \\
acx_4+\left(-y_4c\right)\sqrt{-ab}=\left(aux_1+vaby_1\right)+\left(uy_1-avx_1\right)\sqrt{-ab}.  \label{d}
\end{subnumcases}

{\bf Case 1:} $u\neq\pm c$.

Then $v\neq0$. By \eqref{b} and \eqref{c}, we have $y_2c=uy_2-avx_2$ and $y_3c=uy_3-avx_3$. Then
\begin{equation}\label{e13}
y_2=\frac{av}{u-c}x_2\,\mbox{~and~}\,y_3=\frac{av}{u-c}x_3.
\end{equation}
From \eqref{a} and \eqref{d}, it follows that $-y_1c=uy_4-avx_4$ and $-y_4c=uy_1-avx_1$, i.e.,
$$
\begin{pmatrix} c & u \\ u & c \end{pmatrix}\begin{pmatrix} y_1 \\ y_4 \end{pmatrix}=\begin{pmatrix} avx_4 \\ avx_1 \end{pmatrix}.
$$
Then we have
\begin{equation}\label{e14}
y_1=\frac{av}{c^2-u^2}\left(cx_4-ux_1\right)\,\mbox{~and~}\,y_4=\frac{av}{c^2-u^2}\left(cx_1-ux_4\right).
\end{equation}
Therefore, from \eqref{e13}, \eqref{e14} and $u^2+v^2ab=c^2$, it follows that
\begin{equation}\label{e15}
Y=\frac{av}{c^2-u^2}\begin{pmatrix}cx_4-ux_1 & -\left(u+c\right)x_2 \\ -\left(u+c\right)x_3 & cx_1-ux_4 \end{pmatrix}=\frac{1}{vb}\begin{pmatrix}cx_4-ux_1 & -\left(u+c\right)x_2 \\ -\left(u+c\right)x_3 & cx_1-ux_4 \end{pmatrix}.
\end{equation}
Since $Y\in M_2(\mathbb{Z})$, we have $y_i\in\mathbb{Z},\,i=1,\,2,\,3,\,4$. Then
\begin{subnumcases}
{}vb\mid\left(cx_4-ux_1\right), \label{a1}\\
vb\mid\left(-\left(u+c\right)x_2\right), \label{b1} \\
vb\mid\left(-\left(u+c\right)x_3\right), \label{c1} \\
vb\mid\left(cx_1-ux_4\right).  \label{d1}
\end{subnumcases}
From \eqref{a1} and \eqref{d1}, we obtain $vb\mid\left(-\left(u+c\right)\left(x_1-x_4\right)\right)$. Then there is an integer $s$ such that $-\left(u+c\right)\left(x_1-x_4\right)=vbs$. Since $u^2+v^2ab=c^2$, we have
$$
x_1-x_4=\frac{-vbs}{u+c}=\frac{vbs\left(u-c\right)}{c^2-u^2}=\frac{s\left(u-c\right)}{va}.
$$
Since $x_1-x_4\in\mathbb{Z}$, we obtain $va\mid\left(s\left(u-c\right)\right)$. Let $g=\gcd\left(va,\,u-c\right)$. Then $\frac{va}{g}\mid s$. So there is an integer $k$ such that $s=\frac{va}{g}k$. Thus,
\begin{equation}\label{e18}
x_1-x_4=\frac{s\left(u-c\right)}{va}=\frac{va}{g}k\cdot\frac{\left(u-c\right)}{va}=\frac{u-c}{g}k.
\end{equation}
Likewise, there are integers $t_2$ and $t_3$ such that
\begin{equation}\label{e19}
x_2=\frac{u-c}{g}t_2\,\mbox{~and~}\,x_3=\frac{u-c}{g}t_3.
\end{equation}
From \eqref{e18} and \eqref{e19}, we conclude that
\begin{equation}\label{e16}
X=x_1I+\frac{u-c}{g}\begin{pmatrix} 0 & t_2 \\ t_3 & -k \end{pmatrix}.
\end{equation}
From $-\left(u+c\right)\left(x_1-x_4\right)=vbs$, $s=\frac{va}{g}k$ and \eqref{e15}, it follows that
\begin{equation}\label{e20}
\frac{va}{g}k=s=\frac{-\left(u+c\right)\left(x_1-x_4\right)}{vb}=\frac{cx_4-ux_1-\left(cx_1-ux_4\right)}{vb}=y_1-y_4.
\end{equation}
From $x_2=\frac{u-c}{g}t_2$, $u^2+v^2ab=c^2$ and \eqref{e15}, it follows that
\begin{equation}\label{e21}
y_2=\frac{-\left(u+c\right)x_2}{vb}=\frac{\left(c^2-u^2\right)x_2}{vb\left(u-c\right)}=\frac{va}{u-c}x_2=\frac{va}{u-c}\cdot\frac{u-c}{g}t_2=\frac{va}{g}t_2.
\end{equation}
Similarly, we obtain
\begin{equation}\label{e22}
y_3=\frac{va}{g}t_3.
\end{equation}
From \eqref{e20}, \eqref{e21} and \eqref{e22}, we conclude that
\begin{equation}\label{e17}
Y=y_1I+\frac{va}{g}\begin{pmatrix} 0 & t_2 \\ t_3 & -k \end{pmatrix}.
\end{equation}
By \eqref{e16} and \eqref{e17}, we have
$$
aX^2+bY^2=tI+\frac{2a}{g^2}r\begin{pmatrix} 0 & t_2 \\ t_3 & -k \end{pmatrix},
$$
where $t=ax_1^2+by_1^2+\frac{2act_2t_3}{g^2}\left(c-u\right)$ and $r=\left(gx_1+ck\right)\left(u-c\right)+vbgy_1$.
From the assumption $aX^2+bY^2=cI$, it follows that $t=c$ and $rk=rt_2=rt_3=0$. If $k=t_2=t_3=0$, then $X=x_1I$ and $Y=y_1I$, where $ax_1^2+by_1^2=c$. If $k,\,t_2,\,t_3$ are not all equal to zero, then $r=0$ and $t=c$.

{\bf Case 2:} $u=-c$.

Then $v=0$. From \eqref{a}, \eqref{b} and \eqref{c}, it follows that $x_1=-x_4$, $y_1=y_4$ and $y_2=y_3=0$. Then $X=\begin{pmatrix} x_1 & x_2 \\ x_3 & -x_1 \end{pmatrix}$ and $Y=y_1I$. The assumption $aX^2+bY^2=cI$ implies that $a\left(x_1^2+x_2x_3\right)+by_1^2=c$. Note that the matrices and the condition which we obtain in this case can be obtained by taking $u=-c$ in Case 1. Indeed, let $u=-c$ in Case 1. Then $v=0$ and $g=\gcd\left(va,\,u-c\right)=2|c|$. By Case 1, we obtain
\[
X=\begin{pmatrix} t_1 & -\frac{c}{|c|}t_2 \\ -\frac{c}{|c|}t_3 & -t_1 \end{pmatrix}\,\mbox{~and~}\,Y=t_4I,
\]
where $t_1,\,t_2,\,t_3,\,t_4\in\mathbb{Z}$ satisfy $a\left(t_1^2+t_2t_3\right)+bt_4^2=c$. If $c>0$, then let $t_1=x_1$, $t_2=-x_2$, $t_3=-x_3$ and $t_4=y_1$. Otherwise, let $t_1=x_1$, $t_2=x_2$, $t_3=x_3$ and $t_4=y_1$. Then we can obtain $X=\begin{pmatrix} x_1 & x_2 \\ x_3 & -x_1 \end{pmatrix}$ and $Y=y_1I$, where $a\left(x_1^2+x_2x_3\right)+by_1^2=c$.

{\bf Case 3:} $u=c$.

Then $v=0$. From \eqref{a}, \eqref{b} and \eqref{c}, it follows that $x_1=x_4$, $y_1=-y_4$ and $x_2=x_3=0$. Then $X=x_1I$ and $Y=\begin{pmatrix} y_1 & y_2 \\ y_3 & -y_1 \end{pmatrix}$. The assumption $aX^2+bY^2=cI$ implies that $ax_1^2+b\left(y_1^2+y_2y_3\right)=c$.
\end{prof}

About Theorem \ref{th3}, we have the following equivalent statement.

\begin{theorem}\label{co1}
Let $a,\, b,\, c$ be nonzero integers such that $-ab$ is not a square and $\gcd\left(a,\, b,\, c\right)=1$. Let $X$ and $Y$ be $2\times2$ matrices over $\mathbb{Z}$. Then
$$
aX^2+bY^2=cI,\,XY=YX
$$
if and only if one of the following statements holds.
\begin{enumerate}
\item[\rm(i)]\, $X=t_1I$ and $Y=t_2I$, where $t_1,\,t_2\in\mathbb{Z}$ satisfy $at_1^2+bt_2^2=c$;
\item[\rm(ii)]\, $X=t_1I$ and $Y=\begin{pmatrix} t_4 & t_2 \\ t_3 & -t_4 \end{pmatrix}$, where $t_1,\,t_2,\,t_3,\,t_4\in\mathbb{Z}$ satisfy $at_1^2+b\left(t_4^2+t_2t_3\right)=c$;
\item[\rm(iii)]\,$X=\begin{pmatrix} t_1 & \frac{u-c}{g}t_2 \\ \frac{u-c}{g}t_3 & \frac{ut_1+vbt_4}{c} \end{pmatrix}$ and $Y=\begin{pmatrix} t_4 & \frac{va}{g}t_2 \\ \frac{va}{g}t_3 & \frac{vat_1-ut_4}{c} \end{pmatrix}$, where $t_1,\,t_2,\,t_3,\,t_4,\,u,\,v\in\mathbb{Z},\,u\neq c,\,g=\gcd\left(va,\,u-c\right)$ and the following relations hold:
$$
u^2+v^2ab=c^2,\,at_1^2+bt_4^2+\frac{2act_2t_3}{g^2}\left(c-u\right)=c,\,c\mid\left(ut_1+vbt_4\right),\,c\mid\left(vat_1-ut_4\right).
$$
\end{enumerate}
\end{theorem}

\begin{prof}
We only need to show that (iii) is equivalent to Theorem \ref{th3} (iii).

We now prove sufficiency. By Theorem \ref{th3} (iii), we have $k=\frac{g}{c\left(c-u\right)}\left[\left(u-c\right)t_1+vbt_4\right]$. Then
\[
X=t_1I+\frac{u-c}{g}\begin{pmatrix} 0 & t_2 \\ t_3 & -k \end{pmatrix}=\begin{pmatrix} t_1 & \frac{u-c}{g}t_2 \\ \frac{u-c}{g}t_3 & \frac{ut_1+vbt_4}{c} \end{pmatrix}
\]
and
\[
Y=t_4I+\frac{va}{g}\begin{pmatrix} 0 & t_2 \\ t_3 & -k \end{pmatrix}=\begin{pmatrix} t_4 & \frac{va}{g}t_2 \\ \frac{va}{g}t_3 & \frac{vat_1-ut_4}{c} \end{pmatrix}.
\]
Since $k\in\mathbb{Z}$, we have $c\mid\left(ut_1+vbt_4\right)$ and $c\mid\left(vat_1-ut_4\right)$.

We next prove necessity. By (iii), we have
$$X=\begin{pmatrix} t_1 & \frac{u-c}{g}t_2 \\ \frac{u-c}{g}t_3 & \frac{ut_1+vbt_4}{c} \end{pmatrix}=t_1I+\frac{u-c}{g}\begin{pmatrix} 0 & t_2 \\ t_3 & -k \end{pmatrix}$$
and
$$Y=\begin{pmatrix} t_4 & \frac{va}{g}t_2 \\ \frac{va}{g}t_3 & \frac{vat_1-ut_4}{c} \end{pmatrix}=t_4I+\frac{va}{g}\begin{pmatrix} 0 & t_2 \\ t_3 & -k \end{pmatrix},$$
where $k=\frac{g}{c\left(c-u\right)}\left[\left(u-c\right)t_1+vbt_4\right]$. To complete the proof we need to show that $k\in\mathbb{Z}$ and $\left(gt_1+ck\right)\left(u-c\right)+vbgt_4=0$. Since $c\mid\left(ut_1+vbt_4\right)$ and $c\mid\left(vat_1-ut_4\right)$, we have $\frac{u-c}{g}k,\,\frac{va}{g}k\in\mathbb{Z}$. Then $c\mid\left(\left(u-c\right)t_1+vbt_4\right)$ and $\frac{c-u}{g}\mid\left(\frac{va}{g}\cdot\frac{\left(u-c\right)t_1+vbt_4}{c}\right)$. Since $\gcd\left(\frac{va}{g},\,\frac{c-u}{g}\right)=1$, we obtain $\frac{c-u}{g}\mid\frac{\left(u-c\right)t_1+vbt_4}{c}$. Hence, $k\in\mathbb{Z}$. Moreover, by a direct computation, we have $\left(gt_1+ck\right)\left(u-c\right)+vbgt_4=0$.
\end{prof}

By Theorem \ref{co1}, we can get all solutions of some matrix equations for given nonzero integers $a,\, b,\, c$.

\begin{proposition}
Let $p$ be a prime such that $p\equiv3\pmod 4$. Then all solutions of the matrix equations
\begin{equation}\label{e24}
X^2+Y^2=\pm pI,\,X,\,Y\in M_2(\mathbb{Z})
\end{equation}
are given by the following five parts.
\begin{enumerate}
\item[\rm(i)]\,$X=\begin{pmatrix} t_1 & t_2 \\ t_3 & -t_1 \end{pmatrix}$ and $Y=\begin{pmatrix} s_1 & s_2 \\ s_3 & -s_1 \end{pmatrix}$, where $t_1,\,t_2,\,t_3,\,s_1,\,s_2,\,s_3\in\mathbb{Z}$ satisfy $t_1^2+t_2t_3+s_1^2+s_2s_3=\pm p$;
\item[\rm(ii)]\, $X=t_1I$ and $Y=\begin{pmatrix} t_4 & t_2 \\ t_3 & -t_4 \end{pmatrix}$, where $t_1,\,t_2,\,t_3,\,t_4\in\mathbb{Z}$ satisfy $t_1^2+t_4^2+t_2t_3=\pm p$;
\item[\rm(iii)]\,$X=\begin{pmatrix} t_1 & t_2 \\ t_3 & -t_1 \end{pmatrix}$ and $Y=t_4I$, where $t_1,\,t_2,\,t_3,\,t_4\in\mathbb{Z}$ satisfy $t_1^2+t_4^2+t_2t_3=\pm p$;
\item[\rm(iv)]\,$X=\begin{pmatrix} t_1 & t_2 \\ t_3 & t_4 \end{pmatrix}$ and $Y=\begin{pmatrix} t_4 & -t_2 \\ -t_3 & t_1 \end{pmatrix}$, where $t_1,\,t_2,\,t_3,\,t_4\in\mathbb{Z}$ satisfy $t_1^2+t_4^2+2t_2t_3=\pm p$;
\item[\rm(v)]\,$X=\begin{pmatrix} t_1 & t_2 \\ t_3 & -t_4 \end{pmatrix}$ and $Y=\begin{pmatrix} t_4 & t_2 \\ t_3 & -t_1 \end{pmatrix}$, where $t_1,\,t_2,\,t_3,\,t_4\in\mathbb{Z}$ satisfy $t_1^2+t_4^2+2t_2t_3=\pm p$.
\end{enumerate}
\end{proposition}

\begin{prof}
We can prove this proposition simultaneously for the equations $X^2+Y^2=pI$ and $X^2+Y^2=-pI$, where the upper signs refer to the first equation and the lower signs refer to the second equation. From Proposition \ref{p2} 1), it follows that all non-commutative solutions of equations \eqref{e24} are given by (i). We next find commutative solutions of equations \eqref{e24}. By Theorem \ref{co1}, we only need to consider the following three cases.

{\bf Case 1:} $X=t_1I$ and $Y=t_2I$, where $t_1,\,t_2\in\mathbb{Z}$ satisfy $t_1^2+t_2^2=\pm p$.

Obviously, $t_1^2+t_2^2=-p$ is impossible. Since $p$ is a prime such that $p\equiv3\pmod 4$, it follows that $p$ cannot be represented as a sum of two squares. This means that $t_1^2+t_2^2=p$ is impossible. Therefore, this case is impossible.

{\bf Case 2:} $X=t_1I$ and $Y=\begin{pmatrix} t_4 & t_2 \\ t_3 & -t_4 \end{pmatrix}$, where $t_1,\,t_2,\,t_3,\,t_4\in\mathbb{Z}$ satisfy $t_1^2+t_4^2+t_2t_3=\pm p$.

Then we can get the solutions (ii) in this case.

{\bf Case 3:} $X=\begin{pmatrix} t_1 & \frac{u\mp p}{g}t_2 \\ \frac{u\mp p}{g}t_3 & \frac{ut_1+vt_4}{\pm p} \end{pmatrix}$ and $Y=\begin{pmatrix} t_4 & \frac{v}{g}t_2 \\ \frac{v}{g}t_3 & \frac{vt_1-ut_4}{\pm p} \end{pmatrix}$, where $t_1,\,t_2,\,t_3,\,t_4,\,u,\,v\in\mathbb{Z},\,u\neq \pm p,\,g=\gcd\left(v,\,u\mp p\right)$ and the following relations hold:
$$
u^2+v^2=p^2,\,t_1^2+t_4^2+\frac{2pt_2t_3}{g^2}\left(p\mp u\right)=\pm p,\,p\mid\left(ut_1+vt_4\right),\,p\mid\left(vt_1-ut_4\right).
$$

In this case, we need to solve the Diophantine equation $u^2+v^2=p^2$ in integers $u,\,v$. Note that $p$ is a prime such that $p\equiv3\pmod 4$. Then it follows from $u^2+v^2=p^2$ that $p\mid u$ and $p\mid v$. Thus, we obtain $\left(u,\,v\right)\in\{\left(\mp p,\,0\right),\,\left(0,\,p\right),\,\left(0,\,-p\right)\}$.

We now consider the matrix equation $X^2+Y^2=pI$. Then we have
\[
\left(u,\,v\right)\in\{\left(-p,\,0\right),\,\left(0,\,p\right),\,\left(0,\,-p\right)\}.
\]
For $\left(u,\,v\right)=\left(-p,\,0\right)$, we have $X=\begin{pmatrix} t_1 & t_2 \\ t_3 & -t_1 \end{pmatrix}$ and $Y=t_4I$, where $t_1,\,t_2,\,t_3,\,t_4\in\mathbb{Z}$ satisfy $t_1^2+t_4^2+t_2t_3=p$. For $\left(u,\,v\right)=\left(0,\,p\right)$, we have $X=\begin{pmatrix} t_1 & t_2 \\ t_3 & t_4 \end{pmatrix}$ and $Y=\begin{pmatrix} t_4 & -t_2 \\ -t_3 & t_1 \end{pmatrix}$, where $t_1,\,t_2,\,t_3,\,t_4\in\mathbb{Z}$ satisfy $t_1^2+t_4^2+2t_2t_3=p$. For $\left(u,\,v\right)=\left(0,\,-p\right)$, we have $X=\begin{pmatrix} t_1 & t_2 \\ t_3 & -t_4 \end{pmatrix}$ and $Y=\begin{pmatrix} t_4 & t_2 \\ t_3 & -t_1 \end{pmatrix}$, where $t_1,\,t_2,\,t_3,\,t_4\in\mathbb{Z}$ satisfy $t_1^2+t_4^2+2t_2t_3=p$.

We next consider the matrix equation $X^2+Y^2=-pI$. Then we obtain
\[
\left(u,\,v\right)\in\{\left(p,\,0\right),\,\left(0,\,-p\right),\,\left(0,\,p\right)\}.
\]
For $\left(u,\,v\right)=\left(p,\,0\right)$, we have $X=\begin{pmatrix} t_1 & t_2 \\ t_3 & -t_1 \end{pmatrix}$ and $Y=t_4I$, where $t_1,\,t_2,\,t_3,\,t_4\in\mathbb{Z}$ satisfy $t_1^2+t_4^2+t_2t_3=-p$. For $\left(u,\,v\right)=\left(0,\,-p\right)$, we have $X=\begin{pmatrix} t_1 & t_2 \\ t_3 & t_4 \end{pmatrix}$ and $Y=\begin{pmatrix} t_4 & -t_2 \\ -t_3 & t_1 \end{pmatrix}$, where $t_1,\,t_2,\,t_3,\,t_4\in\mathbb{Z}$ satisfy $t_1^2+t_4^2+2t_2t_3=-p$. For $\left(u,\,v\right)=\left(0,\,p\right)$, we have $X=\begin{pmatrix} t_1 & t_2 \\ t_3 & -t_4 \end{pmatrix}$ and $Y=\begin{pmatrix} t_4 & t_2 \\ t_3 & -t_1 \end{pmatrix}$, where $t_1,\,t_2,\,t_3,\,t_4\in\mathbb{Z}$ satisfy $t_1^2+t_4^2+2t_2t_3=-p$.
\end{prof}

\begin{proposition}
Let $p$ be a prime such that $p\equiv5\,\text{or}\,7\pmod 8$. Then all solutions of the matrix equations
\begin{equation}\label{e25}
X^2+2Y^2=\pm pI,\,X,\,Y\in M_2(\mathbb{Z})
\end{equation}
are given by the following three parts.
\begin{enumerate}
\item[\rm(i)]\,$X=\begin{pmatrix} t_1 & t_2 \\ t_3 & -t_1 \end{pmatrix}$ and $Y=\begin{pmatrix} s_1 & s_2 \\ s_3 & -s_1 \end{pmatrix}$, where $t_1,\,t_2,\,t_3,\,s_1,\,s_2,\,s_3\in\mathbb{Z}$ satisfy $t_1^2+t_2t_3+2\left(s_1^2+s_2s_3\right)=\pm p$;
\item[\rm(ii)]\, $X=t_1I$ and $Y=\begin{pmatrix} t_4 & t_2 \\ t_3 & -t_4 \end{pmatrix}$, where $t_1,\,t_2,\,t_3,\,t_4\in\mathbb{Z}$ satisfy $t_1^2+2\left(t_4^2+t_2t_3\right)=\pm p$;
\item[\rm(iii)]\, $X=\begin{pmatrix} t_1 & t_2 \\ t_3 & -t_1 \end{pmatrix}$ and $Y=t_4I$, where $t_1,\,t_2,\,t_3,\,t_4\in\mathbb{Z}$ satisfy $t_1^2+t_2t_3+2t_4^2=\pm p$.

\end{enumerate}
\end{proposition}

\begin{prof}
We can prove this proposition simultaneously for the equations $X^2+2Y^2=pI$ and $X^2+2Y^2=-pI$, where the upper signs refer to the first equation and the lower signs refer to the second equation. From Proposition \ref{p2} 1), it follows that all non-commutative solutions of equations \eqref{e25} are given by (i). We next find commutative solutions of equations \eqref{e25}. By Theorem \ref{co1}, we only need to consider the following three cases.

{\bf Case 1:} $X=t_1I$ and $Y=t_2I$, where $t_1,\,t_2\in\mathbb{Z}$ satisfy $t_1^2+2t_2^2=\pm p$.

Obviously, $t_1^2+2t_2^2=-p$ is impossible. We claim that $t_1^2+2t_2^2=p$ is also impossible. Indeed, if $t_1^2+2t_2^2=p$, then $\gcd\left(p,\,t_1t_2\right)=1$. Otherwise, we have $p\mid t_1$ and $p\mid t_2$. This is impossible by $t_1^2+2t_2^2=p$. So $\gcd\left(p,\,t_1t_2\right)=1$. Then there is an integer $t_2'$ such that $t_2t_2'\equiv1\pmod p$. From $t_1^2+2t_2^2=p$, it follows that $\left(t_1t_2'\right)^2+2\left(t_2t_2'\right)^2=p\left(t_2'\right)^2$. Then $\left(t_1t_2'\right)^2\equiv-2\pmod p$. This means that $\left(\frac{-2}{p}\right)=1$, where $\left(\frac{\cdot}{p}\right)$ is the Legendre symbol. However, for $p\equiv5\,\text{or}\,7\pmod 8$, we have $\left(\frac{-2}{p}\right)=-1$, a contradiction.

{\bf Case 2:} $X=t_1I$ and $Y=\begin{pmatrix} t_4 & t_2 \\ t_3 & -t_4 \end{pmatrix}$, where $t_1,\,t_2,\,t_3,\,t_4\in\mathbb{Z}$ satisfy $t_1^2+2\left(t_4^2+t_2t_3\right)=\pm p$.

Then we can get the solutions (ii) in this case.

{\bf Case 3:} $X=\begin{pmatrix} t_1 & \frac{u\mp p}{g}t_2 \\ \frac{u\mp p}{g}t_3 & \frac{ut_1+2vt_4}{\pm p} \end{pmatrix}$ and $Y=\begin{pmatrix} t_4 & \frac{v}{g}t_2 \\ \frac{v}{g}t_3 & \frac{vt_1-ut_4}{\pm p} \end{pmatrix}$, where $t_1,\,t_2,\,t_3,\,t_4,\,u,\,v\in\mathbb{Z},\,u\neq\pm p,\,g=\gcd\left(v,\,u\mp p\right)$ and the following relations hold:
$$
u^2+2v^2=p^2,\,t_1^2+2t_4^2+\frac{2pt_2t_3}{g^2}\left(p\mp u\right)=\pm p,\,p\mid\left(ut_1+2vt_4\right),\,p\mid\left(vt_1-ut_4\right).
$$

In this case, we need to solve the Diophantine equation $u^2+2v^2=p^2$ in integers $u,\,v$. From $u^2+2v^2=p^2$, it follows that $2\mid\left(p-u\right)\left(p+u\right)$. Then $p$ and $u$ have the same parity. Moreover, $2\mid v$. By $u^2+2v^2=p^2$, we have
\begin{equation}\label{e26}
2\left(\frac{v}{2}\right)^2=\frac{p-u}{2}\cdot\frac{p+u}{2}.
\end{equation}
If $v=0$, then we obtain $\left(u,\,v\right)=\left(\mp p,\,0\right)$. Let us now assume that $v\neq0$. Let $g=\gcd\left(\frac{p-u}{2},\,\frac{p+u}{2}\right)$. Then $g\mid p$. If $g=1$, then it follows from \eqref{e26} that there are integers $y_1$ and $y_2$ such that
\[
\frac{p-u}{2}=2y_1^2,\,\frac{p+u}{2}=y_2^2\,\mbox{~or~}\,\frac{p-u}{2}=y_1^2,\,\frac{p+u}{2}=2y_2^2.
\]
Then we have $p=2y_1^2+y_2^2$ or $p=y_1^2+2y_2^2$, which is impossible by the argument of Case 1. Therefore, $g=p$. Then $p\mid u$, so $p\mid v$. From $u^2+2v^2=p^2$, it follows that $\left(u,\,v\right)=\left(\mp p,\,0\right)$. Then we have $X=\begin{pmatrix} t_1 & t_2 \\ t_3 & -t_1 \end{pmatrix}$ and $Y=t_4I$, where $t_1,\,t_2,\,t_3,\,t_4\in\mathbb{Z}$ satisfy $t_1^2+t_2t_3+2t_4^2=\pm p$.
\end{prof}

Let $c=\pm1$ in Theorem \ref{co1}. Then we have the following corollary.
\begin{corollary}\label{co3}
Let $a,\, b$ be nonzero integers such that $-ab$ is not a square and let $X,\,Y\in M_2(\mathbb{Z})$. Then
$$
aX^2+bY^2=\pm I,\,XY=YX
$$
if and only if one of the following statements holds.
\begin{enumerate}
\item[\rm(i)]\, $X=t_1I$ and $Y=t_2I$, where $t_1,\,t_2\in\mathbb{Z}$ satisfy $at_1^2+bt_2^2=\pm1$;
\item[\rm(ii)]\, $X=t_1I$ and $Y=\begin{pmatrix} t_4 & t_2 \\ t_3 & -t_4 \end{pmatrix}$, where $t_1,\,t_2,\,t_3,\,t_4\in\mathbb{Z}$ satisfy $at_1^2+b\left(t_4^2+t_2t_3\right)=\pm1$;
\item[\rm(iii)]\,$X=\begin{pmatrix} t_1 & \frac{u\mp1}{g}t_2 \\ \frac{u\mp1}{g}t_3 & \pm\left(ut_1+vbt_4\right) \end{pmatrix}$ and $Y=\begin{pmatrix} t_4 & \frac{va}{g}t_2 \\ \frac{va}{g}t_3 & \pm\left(vat_1-ut_4\right) \end{pmatrix}$, where $t_1,\,t_2,\,t_3,\,t_4,\,u,\,v\in\mathbb{Z},\,u\neq \pm1,\,g=\gcd\left(va,\,u\mp1\right)$ and the following relations hold:
$$
u^2+v^2ab=1,\,at_1^2+bt_4^2+\frac{2at_2t_3}{g^2}\left(1\mp u\right)=\pm1.
$$
\end{enumerate}
\end{corollary}

\begin{remark}
\rm Let $d$ be a square-free integer and let $a=1,\,b=-d$. Then Corollary \ref{co3} becomes \cite[Theorem 2.1]{Cohen}.
\end{remark}

By Corollary \ref{co3}, we have the following proposition.

\begin{proposition}\label{p4}
Let $b$ be a nonzero integer such that $-b$ is not a square. Suppose that $b$ has a prime divisor of the form $4k+3$. Then all solutions of the matrix equation
\begin{equation}\label{e27}
X^2+bY^2=-I,\,X,\,Y\in M_2(\mathbb{Z})
\end{equation}
are given by the following two parts.
\begin{enumerate}
\item[\rm(i)]\,$X=\begin{pmatrix} t_1 & t_2 \\ t_3 & -t_1 \end{pmatrix}$ and $Y=\begin{pmatrix} s_1 & s_2 \\ s_3 & -s_1 \end{pmatrix}$, where $t_1,\,t_2,\,t_3,\,s_1,\,s_2,\,s_3\in\mathbb{Z}$ satisfy $t_1^2+t_2t_3+b\left(s_1^2+s_2s_3\right)=-1$;
\item[\rm(ii)]\,For $b>0$, let $X=\begin{pmatrix} t_1 & t_2 \\ t_3 & -t_1 \end{pmatrix}$ and $Y=t_4I$, where $t_1,\,t_2,\,t_3,\,t_4\in\mathbb{Z}$ satisfy $t_1^2+t_2t_3+bt_4^2=-1$. Otherwise, let $X=\begin{pmatrix} t_1 & \frac{u+1}{g}t_2 \\ \frac{u+1}{g}t_3 & -\left(ut_1+vbt_4\right) \end{pmatrix}$ and $Y=\begin{pmatrix} t_4 & \frac{v}{g}t_2 \\ \frac{v}{g}t_3 & ut_4-vt_1 \end{pmatrix}$, where $t_1,\,t_2,\,t_3,\,t_4,\,u,\,v\in\mathbb{Z},\,u\neq-1,\,g=\gcd\left(v,\,u+1\right)$ and the following relations hold:
$$
u^2+v^2b=1,\,t_1^2+bt_4^2+\frac{2t_2t_3}{g^2}\left(1+u\right)=-1.
$$
\end{enumerate}
\end{proposition}
\begin{prof}
From Proposition \ref{p2} 1), it follows that all non-commutative solutions of equation \eqref{e27} are given by (i). We next find commutative solutions of equation \eqref{e27}. Suppose that $p\equiv3\pmod4$ is the prime divisor of $b$. Then $\left(\frac{-1}{p}\right)=-1$, where $\left(\frac{\cdot}{p}\right)$ is the Legendre symbol. By Corollary \ref{co3}, we only need to consider the following three cases.

{\bf Case 1:} $X=t_1I$ and $Y=t_2I$, where $t_1,\,t_2\in\mathbb{Z}$ satisfy $t_1^2+bt_2^2=-1$.

We claim that $t_1^2+bt_2^2=-1$ is impossible. Indeed, if $t_1^2+bt_2^2=-1$, then $t_1^2\equiv-1\pmod p$. This means that $\left(\frac{-1}{p}\right)=1$, a contradiction.

{\bf Case 2:} $X=t_1I$ and $Y=\begin{pmatrix} t_4 & t_2 \\ t_3 & -t_4 \end{pmatrix}$, where $t_1,\,t_2,\,t_3,\,t_4\in\mathbb{Z}$ satisfy $t_1^2+b\left(t_4^2+t_2t_3\right)=-1$.

However, this case is also impossible, where the reason is similar to Case 1.

{\bf Case 3:} $X=\begin{pmatrix} t_1 & \frac{u+1}{g}t_2 \\ \frac{u+1}{g}t_3 & -\left(ut_1+vbt_4\right) \end{pmatrix}$ and $Y=\begin{pmatrix} t_4 & \frac{v}{g}t_2 \\ \frac{v}{g}t_3 & ut_4-vt_1 \end{pmatrix}$, where $t_1,\,t_2,\,t_3,\,t_4,\,u,\,v\in\mathbb{Z},\,u\neq-1,\,g=\gcd\left(v,\,u+1\right)$ and the following relations hold:
$$
u^2+v^2b=1,\,t_1^2+bt_4^2+\frac{2t_2t_3}{g^2}\left(1+u\right)=-1.
$$

If $b>0$, then it follows from $u^2+v^2b=1$ that $u=1$ and $v=0$. Then we have $X=\begin{pmatrix} t_1 & t_2 \\ t_3 & -t_1 \end{pmatrix}$ and $Y=t_4I$, where $t_1,\,t_2,\,t_3,\,t_4\in\mathbb{Z}$ satisfy $t_1^2+t_2t_3+bt_4^2=-1$. If $b<0$, then $u^2+v^2b=1$ is the Pell's equation. We know that it has infinitely many solutions in integers $u$ and $v$.
\end{prof}

\begin{example}
\rm
Let $b=-3$ in Proposition \ref{p4}. Then all solutions of the matrix equation
\begin{equation}\label{e28}
X^2-3Y^2=-I,\,X,\,Y\in M_2(\mathbb{Z})
\end{equation}
are given by the following two parts.
\begin{enumerate}
\item[\rm(i)]\,$X=\begin{pmatrix} t_1 & t_2 \\ t_3 & -t_1 \end{pmatrix}$ and $Y=\begin{pmatrix} s_1 & s_2 \\ s_3 & -s_1 \end{pmatrix}$, where $t_1,\,t_2,\,t_3,\,s_1,\,s_2,\,s_3\in\mathbb{Z}$ satisfy $t_1^2+t_2t_3-3\left(s_1^2+s_2s_3\right)=-1$;
\item[\rm(ii)]\,$X=\begin{pmatrix} t_1 & \frac{u+1}{g}t_2 \\ \frac{u+1}{g}t_3 & 3vt_4-ut_1 \end{pmatrix}$ and $Y=\begin{pmatrix} t_4 & \frac{v}{g}t_2 \\ \frac{v}{g}t_3 & ut_4-vt_1 \end{pmatrix}$, where $t_1,\,t_2,\,t_3,\,t_4,\,u,\,v\in\mathbb{Z},\,u\neq-1,\,g=\gcd\left(v,\,u+1\right)$ and the following relations hold:
$$
u^2-3v^2=1,\,t_1^2-3t_4^2+\frac{2t_2t_3}{g^2}\left(1+u\right)=-1.
$$
\end{enumerate}

For example, let $u=7$ and $v=4$ in (ii). Then we can get a family of solutions to equation \eqref{e28}:
\[
\begin{pmatrix} t_1 & 2t_2 \\ 2t_3 & 12t_4-7t_1 \end{pmatrix}^2-3\begin{pmatrix} t_4 & t_2 \\ t_3 & 7t_4-4t_1 \end{pmatrix}^2=-I,
\]
where $t_1,\,t_2,\,t_3,\,t_4\in\mathbb{Z}$ satisfy $t_1^2-3t_4^2+t_2t_3=-1$.
\end{example}

Let $c=\pm2$ in Theorem \ref{co1}. Then we have the following corollary.

\begin{corollary}\label{co2}
Let $a,\, b$ be nonzero integers such that $-ab$ is not a square and let $X,\,Y\in M_2(\mathbb{Z})$. Then
$$
aX^2+bY^2=\pm2I,\,XY=YX
$$
if and only if one of the following statements holds.
\begin{enumerate}
\item[\rm(i)]\, $X=t_1I$ and $Y=t_2I$, where $t_1,\,t_2\in\mathbb{Z}$ satisfy $at_1^2+bt_2^2=\pm2$;
\item[\rm(ii)]\, $X=t_1I$ and $Y=\begin{pmatrix} t_4 & t_2 \\ t_3 & -t_4 \end{pmatrix}$, where $t_1,\,t_2,\,t_3,\,t_4\in\mathbb{Z}$ satisfy $at_1^2+b\left(t_4^2+t_2t_3\right)=\pm2$;
\item[\rm(iii)]\,$X=\begin{pmatrix} t_1 & \frac{u\mp2}{g}t_2 \\ \frac{u\mp2}{g}t_3 & \frac{ut_1+vbt_4}{\pm2} \end{pmatrix}$ and $Y=\begin{pmatrix} t_4 & \frac{va}{g}t_2 \\ \frac{va}{g}t_3 & \frac{vat_1-ut_4}{\pm2} \end{pmatrix}$, where $t_1,\,t_2,\,t_3,\,t_4,\,u,\,v\in\mathbb{Z},\,u\neq\pm2,\,g=\gcd\left(va,\,u\mp2\right)$ and the following relations hold:
$$
u^2+v^2ab=4,\,at_1^2+bt_4^2+\frac{4at_2t_3}{g^2}\left(2\mp u\right)=\pm2.
$$
\end{enumerate}
\end{corollary}

\begin{prof}
By Theorem \ref{co1}, we only need to show that $2\mid\left(ut_1+vbt_4\right)$ and $2\mid\left(vat_1-ut_4\right)$ in (iii). To see that $g=\gcd\left(va,\,u\mp2\right)$, the proof is divided into the following two cases, depending on whether $g$ is even or not.

{\bf Case 1:} $2\mid g$.

Then $2\mid va$ and $2\mid u$. Therefore, $2\mid\left(vat_1-ut_4\right)$. If $2\mid v$, then $2\mid\left(ut_1+vbt_4\right)$. Otherwise, by $u^2+v^2ab=4$, we have $4\mid ab$. If $2\mid b$, then $2\mid\left(ut_1+vbt_4\right)$. If $2\nmid b$, then $2\mid a$. From $at_1^2+bt_4^2+\frac{4at_2t_3}{g^2}\left(2\mp u\right)=\pm2$, it follows that $vat_1^2+vbt_4^2+\frac{4vat_2t_3}{g^2}\left(2\mp u\right)=\pm2v$. Then $2\mid t_4$, so $2\mid\left(ut_1+vbt_4\right)$.

{\bf Case 2:} $2\nmid g$.

From $t_1,\,t_2,\,t_3,\,t_4,\,u,\,v\in\mathbb{Z}$ and $at_1^2+bt_4^2+\frac{4at_2t_3}{g^2}\left(2\mp u\right)=\pm2$, it follows that $g^2\mid4at_2t_3\left(2\mp u\right)$. Since $2\nmid g$, we obtain $g^2\mid at_2t_3\left(2\mp u\right)$. Then from $u^2+v^2ab=4$ and $at_1^2+bt_4^2+\frac{4at_2t_3}{g^2}\left(2\mp u\right)=\pm2$, it follows that
\begin{equation}\label{e23}
u^2+v^2ab\equiv0\pmod4\,\mbox{~and~}\,at_1^2+bt_4^2\equiv0\pmod2.
\end{equation}
If $2\mid u$, then it follows from \eqref{e23} and $2\nmid g$ that $2\nmid va$, $2\mid b$ and $2\mid t_1$. In this case, we have $2\mid\left(ut_1+vbt_4\right)$ and $2\mid\left(vat_1-ut_4\right)$. If $2\nmid u$, then it follows from \eqref{e23} that $2\nmid v$ and $2\nmid ab$. Then $t_1$ and $t_4$ have the same parity. Therefore, in this case, we have
$$
ut_1+vbt_4\equiv t_1+t_4\equiv0\pmod2\,\mbox{~and~}\,vat_1-ut_4\equiv t_1+t_4\equiv0\pmod2.
$$
\end{prof}

By Corollary \ref{co2}, we have the following two propositions.
\begin{proposition}\label{p5}
Let $b\neq3$ be a nonzero integer such that $-b$ is not a square. Suppose that $b$ has a prime divisor of the form $8k+3$ or $8k+5$. Then all solutions of the matrix equation
\begin{equation}\label{e29}
X^2+bY^2=2I,\,X,\,Y\in M_2(\mathbb{Z})
\end{equation}
are given by the following two parts.
\begin{enumerate}
\item[\rm(i)]\,$X=\begin{pmatrix} t_1 & t_2 \\ t_3 & -t_1 \end{pmatrix}$ and $Y=\begin{pmatrix} s_1 & s_2 \\ s_3 & -s_1 \end{pmatrix}$, where $t_1,\,t_2,\,t_3,\,s_1,\,s_2,\,s_3\in\mathbb{Z}$ satisfy $t_1^2+t_2t_3+b\left(s_1^2+s_2s_3\right)=2$;
\item[\rm(ii)]\,For $b>0$, let $X=\begin{pmatrix} t_1 & t_2 \\ t_3 & -t_1 \end{pmatrix}$ and $Y=t_4I$, where $t_1,\,t_2,\,t_3,\,t_4\in\mathbb{Z}$ satisfy $t_1^2+t_2t_3+bt_4^2=2$. Otherwise, let $X=\begin{pmatrix} t_1 & \frac{u-2}{g}t_2 \\ \frac{u-2}{g}t_3 & \frac{ut_1+vbt_4}{2} \end{pmatrix}$ and $Y=\begin{pmatrix} t_4 & \frac{v}{g}t_2 \\ \frac{v}{g}t_3 & \frac{vt_1-ut_4}{2} \end{pmatrix}$, where $t_1,\,t_2,\,t_3,\,t_4,\,u,\,v\in\mathbb{Z},\,u\neq2,\,g=\gcd\left(v,\,u-2\right)$ and the following relations hold:
$$
u^2+v^2b=4,\,t_1^2+bt_4^2+\frac{4t_2t_3}{g^2}\left(2-u\right)=2.
$$
\end{enumerate}
\end{proposition}

\begin{prof}
From Proposition \ref{p2} 1), it follows that all non-commutative solutions of equation \eqref{e29} are given by (i). We next find commutative solutions of equation \eqref{e29}. Suppose that $p\equiv3\,\text{or}\,5\pmod8$ is the prime divisor of $b$. Then $\left(\frac{2}{p}\right)=-1$, where $\left(\frac{\cdot}{p}\right)$ is the Legendre symbol. By Corollary \ref{co2}, we only need to consider the following three cases.

{\bf Case 1:} $X=t_1I$ and $Y=t_2I$, where $t_1,\,t_2\in\mathbb{Z}$ satisfy $t_1^2+bt_2^2=2$.

We claim that $t_1^2+bt_2^2=2$ is impossible. Indeed, if $t_1^2+bt_2^2=2$, then $t_1^2\equiv2\pmod p$. This means that $\left(\frac{2}{p}\right)=1$, a contradiction.

{\bf Case 2:} $X=t_1I$ and $Y=\begin{pmatrix} t_4 & t_2 \\ t_3 & -t_4 \end{pmatrix}$, where $t_1,\,t_2,\,t_3,\,t_4\in\mathbb{Z}$ satisfy $t_1^2+b\left(t_4^2+t_2t_3\right)=2$.

However, this case is also impossible, where the reason is similar to Case 1.

{\bf Case 3:} $X=\begin{pmatrix} t_1 & \frac{u-2}{g}t_2 \\ \frac{u-2}{g}t_3 & \frac{ut_1+vbt_4}{2} \end{pmatrix}$ and $Y=\begin{pmatrix} t_4 & \frac{v}{g}t_2 \\ \frac{v}{g}t_3 & \frac{vt_1-ut_4}{2} \end{pmatrix}$, where $t_1,\,t_2,\,t_3,\,t_4,\,u,\,v\in\mathbb{Z},\,u\neq2,\,g=\gcd\left(v,\,u-2\right)$ and the following relations hold:
$$
u^2+v^2b=4,\,t_1^2+bt_4^2+\frac{4t_2t_3}{g^2}\left(2-u\right)=2.
$$

If $b>0$, then it follows from $u^2+v^2b=4$ that $u=-2$ and $v=0$. Then we have $X=\begin{pmatrix} t_1 & t_2 \\ t_3 & -t_1 \end{pmatrix}$ and $Y=t_4I$, where $t_1,\,t_2,\,t_3,\,t_4\in\mathbb{Z}$ satisfy $t_1^2+t_2t_3+bt_4^2=2$. If $b<0$, then $u^2+v^2b=4$ is the Pell's equation. We know that it has infinitely many solutions in integers $u$ and $v$.
\end{prof}

\begin{example}
\rm
Let $b=-5$ in Proposition \ref{p5}. Then all solutions of the matrix equation
\begin{equation}\label{e30}
X^2-5Y^2=2I,\,X,\,Y\in M_2(\mathbb{Z})
\end{equation}
are given by the following two parts.
\begin{enumerate}
\item[\rm(i)]\,$X=\begin{pmatrix} t_1 & t_2 \\ t_3 & -t_1 \end{pmatrix}$ and $Y=\begin{pmatrix} s_1 & s_2 \\ s_3 & -s_1 \end{pmatrix}$, where $t_1,\,t_2,\,t_3,\,s_1,\,s_2,\,s_3\in\mathbb{Z}$ satisfy $t_1^2+t_2t_3-5\left(s_1^2+s_2s_3\right)=2$;
\item[\rm(ii)]\,$X=\begin{pmatrix} t_1 & \frac{u-2}{g}t_2 \\ \frac{u-2}{g}t_3 & \frac{ut_1-5vt_4}{2} \end{pmatrix}$ and $Y=\begin{pmatrix} t_4 & \frac{v}{g}t_2 \\ \frac{v}{g}t_3 & \frac{vt_1-ut_4}{2} \end{pmatrix}$, where $t_1,\,t_2,\,t_3,\,t_4,\,u,\,v\in\mathbb{Z},\,u\neq2,\,g=\gcd\left(v,\,u-2\right)$ and the following relations hold:
$$
u^2-5v^2=4,\,t_1^2-5t_4^2+\frac{4t_2t_3}{g^2}\left(2-u\right)=2.
$$
\end{enumerate}

For example, let $u=18$ and $v=8$ in (ii). Then we can get a family of solutions to equation \eqref{e30}:
\[
\begin{pmatrix} t_1 & 2t_2 \\ 2t_3 & 9t_1-20t_4 \end{pmatrix}^2-5\begin{pmatrix} t_4 & t_2 \\ t_3 & 4t_1-9t_4 \end{pmatrix}^2=2I,
\]
where $t_1,\,t_2,\,t_3,\,t_4\in\mathbb{Z}$ satisfy $t_1^2-5t_4^2-t_2t_3=2$.
\end{example}

\begin{proposition}\label{p6}
Let $b$ be a nonzero integer such that $-b$ is not a square. Suppose that $b$ has a prime divisor of the form $8k+5$ or $8k+7$. Then all solutions of the matrix equation
\begin{equation}\label{e31}
X^2+bY^2=-2I,\,X,\,Y\in M_2(\mathbb{Z})
\end{equation}
are given by the following two parts.
\begin{enumerate}
\item[\rm(i)]\,$X=\begin{pmatrix} t_1 & t_2 \\ t_3 & -t_1 \end{pmatrix}$ and $Y=\begin{pmatrix} s_1 & s_2 \\ s_3 & -s_1 \end{pmatrix}$, where $t_1,\,t_2,\,t_3,\,s_1,\,s_2,\,s_3\in\mathbb{Z}$ satisfy $t_1^2+t_2t_3+b\left(s_1^2+s_2s_3\right)=-2$;
\item[\rm(ii)]\,For $b>0$, let $X=\begin{pmatrix} t_1 & t_2 \\ t_3 & -t_1 \end{pmatrix}$ and $Y=t_4I$, where $t_1,\,t_2,\,t_3,\,t_4\in\mathbb{Z}$ satisfy $t_1^2+t_2t_3+bt_4^2=-2$. Otherwise, let $X=\begin{pmatrix} t_1 & \frac{u+2}{g}t_2 \\ \frac{u+2}{g}t_3 & \frac{ut_1+vbt_4}{-2} \end{pmatrix}$ and $Y=\begin{pmatrix} t_4 & \frac{v}{g}t_2 \\ \frac{v}{g}t_3 & \frac{vt_1-ut_4}{-2} \end{pmatrix}$, where $t_1,\,t_2,\,t_3,\,t_4,\,u,\,v\in\mathbb{Z},\,u\neq-2,\,g=\gcd\left(v,\,u+2\right)$ and the following relations hold:
$$
u^2+v^2b=4,\,t_1^2+bt_4^2+\frac{4t_2t_3}{g^2}\left(2+u\right)=-2.
$$
\end{enumerate}
\end{proposition}

\begin{prof}
From Proposition \ref{p2} 1), it follows that all non-commutative solutions of equation \eqref{e31} are given by (i). We next find commutative solutions of equation \eqref{e31}. Suppose that $p\equiv5\,\text{or}\,7\pmod8$ is the prime divisor of $b$. Then $\left(\frac{-2}{p}\right)=-1$, where $\left(\frac{\cdot}{p}\right)$ is the Legendre symbol. By Corollary \ref{co2}, we only need to consider the following three cases.

{\bf Case 1:} $X=t_1I$ and $Y=t_2I$, where $t_1,\,t_2\in\mathbb{Z}$ satisfy $t_1^2+bt_2^2=-2$.

We claim that $t_1^2+bt_2^2=-2$ is impossible. Indeed, if $t_1^2+bt_2^2=-2$, then $t_1^2\equiv-2\pmod p$. This means that $\left(\frac{-2}{p}\right)=1$, a contradiction.

{\bf Case 2:} $X=t_1I$ and $Y=\begin{pmatrix} t_4 & t_2 \\ t_3 & -t_4 \end{pmatrix}$, where $t_1,\,t_2,\,t_3,\,t_4\in\mathbb{Z}$ satisfy $t_1^2+b\left(t_4^2+t_2t_3\right)=-2$.

However, this case is also impossible, where the reason is similar to Case 1.

{\bf Case 3:} $X=\begin{pmatrix} t_1 & \frac{u+2}{g}t_2 \\ \frac{u+2}{g}t_3 & \frac{ut_1+vbt_4}{-2} \end{pmatrix}$ and $Y=\begin{pmatrix} t_4 & \frac{v}{g}t_2 \\ \frac{v}{g}t_3 & \frac{vt_1-ut_4}{-2} \end{pmatrix}$, where $t_1,\,t_2,\,t_3,\,t_4,\,u,\,v\in\mathbb{Z},\,u\neq-2,\,g=\gcd\left(v,\,u+2\right)$ and the following relations hold:
$$
u^2+v^2b=4,\,t_1^2+bt_4^2+\frac{4t_2t_3}{g^2}\left(2+u\right)=-2.
$$

If $b>0$, then it follows from $u^2+v^2b=4$ that $u=2$ and $v=0$. Then we have $X=\begin{pmatrix} t_1 & t_2 \\ t_3 & -t_1 \end{pmatrix}$ and $Y=t_4I$, where $t_1,\,t_2,\,t_3,\,t_4\in\mathbb{Z}$ satisfy $t_1^2+t_2t_3+bt_4^2=-2$. If $b<0$, then $u^2+v^2b=4$ is the Pell's equation. We know that it has infinitely many solutions in integers $u$ and $v$.
\end{prof}

\begin{example}
\rm
Let $b=-5$ in Proposition \ref{p6}. Then all solutions of the matrix equation
\begin{equation}\label{e32}
X^2-5Y^2=-2I,\,X,\,Y\in M_2(\mathbb{Z})
\end{equation}
are given by the following two parts.
\begin{enumerate}
\item[\rm(i)]\,$X=\begin{pmatrix} t_1 & t_2 \\ t_3 & -t_1 \end{pmatrix}$ and $Y=\begin{pmatrix} s_1 & s_2 \\ s_3 & -s_1 \end{pmatrix}$, where $t_1,\,t_2,\,t_3,\,s_1,\,s_2,\,s_3\in\mathbb{Z}$ satisfy $t_1^2+t_2t_3-5\left(s_1^2+s_2s_3\right)=-2$;
\item[\rm(ii)]\,$X=\begin{pmatrix} t_1 & \frac{u+2}{g}t_2 \\ \frac{u+2}{g}t_3 & \frac{5vt_4-ut_1}{2} \end{pmatrix}$ and $Y=\begin{pmatrix} t_4 & \frac{v}{g}t_2 \\ \frac{v}{g}t_3 & \frac{vt_1-ut_4}{-2} \end{pmatrix}$, where $t_1,\,t_2,\,t_3,\,t_4,\,u,\,v\in\mathbb{Z},\,u\neq-2,\,g=\gcd\left(v,\,u+2\right)$ and the following relations hold:
$$
u^2-5v^2=4,\,t_1^2-5t_4^2+\frac{4t_2t_3}{g^2}\left(2+u\right)=-2.
$$
\end{enumerate}

For example, let $u=-18$ and $v=8$ in (ii). Then we can get a family of solutions to equation \eqref{e32}:
\[
\begin{pmatrix} t_1 & -2t_2 \\ -2t_3 & 20t_4+9t_1 \end{pmatrix}^2-5\begin{pmatrix} t_4 & t_2 \\ t_3 & -\left(4t_1+9t_4\right) \end{pmatrix}^2=-2I,
\]
where $t_1,\,t_2,\,t_3,\,t_4\in\mathbb{Z}$ satisfy $t_1^2-5t_4^2-t_2t_3=-2$.
\end{example}

\end{document}